\documentclass[journal,twoside,web]{ieeecolor}
\usepackage{generic}
\usepackage{cite}
\usepackage{amsmath,amssymb,amsfonts}
\usepackage{algorithmic}
\usepackage[table,xcdraw]{xcolor}
\usepackage{colortbl}
\definecolor{lightgray}{gray}{0.85}
\usepackage{graphicx}
\usepackage{algorithm,algorithmic}
\usepackage{hyperref}
\hypersetup{hidelinks=true}
\usepackage{textcomp}
\usepackage{cite}
\usepackage{amsfonts,amssymb}
\usepackage{bm}
\usepackage{amsmath}
\usepackage{multirow}
\usepackage{algorithmic}
\usepackage{graphicx}
\usepackage{subfigure}
\usepackage{float}
\usepackage{epstopdf}
\usepackage{textcomp}
\usepackage{graphicx} 
\usepackage{epstopdf}
\usepackage{array}
\usepackage[thmmarks,amsmath]{ntheorem}
\newtheorem{mythm}{Theorem}
\newtheorem{mydef}{Definition}
\newtheorem{assumption}{Assumption}
\newtheorem{lemma}{Lemma}
\newtheorem{remark}{Remark}
\newtheorem{proposition}{Proposition}

\newtheorem{corollary}{Corollary}

\usepackage{stfloats}

\usepackage[title]{appendix}

\usepackage[textwidth=1.2cm,textsize=tiny]{todonotes}
\reversemarginpar

\usepackage{array}
\usepackage{enumerate}
\usepackage{booktabs}
\usepackage{algorithm}
\usepackage{algorithmic}
\usepackage{setspace}
\usepackage{float}
\usepackage{cases}

\usepackage{mathrsfs}
\usepackage{hyperref}
\usepackage{textcomp}
\def\BibTeX{{\rm B\kern-.05em{\sc i\kern-.025em b}\kern-.08em
    T\kern-.1667em\lower.7ex\hbox{E}\kern-.125emX}}
\begin{document}
\title{When Does Selfishness Align with Team Goals? A Structural Analysis of Equilibrium and Optimality}
\author{Gehui Xu, Thomas Parisini, \IEEEmembership{Fellow, IEEE},  and Andreas A. Malikopoulos, \IEEEmembership{Senior Member, IEEE}
\thanks{*This research was supported in part by NSF under Grants CNS-2401007, CMMI-2348381, IIS-2415478, in part by MathWorks, and in part by the Horizon Europe Marie Sklodowska-Curie Action under Grant  101204086.
}
\thanks{G. Xu is with the Department of Electrical and Electronic Engineering, Imperial College London, UK. E-mail: {\tt\small g.xu@imperial.ac.uk}}
\thanks{Thomas Parisini is with  the Dept. of Electrical and Electronic Engineering,
	Imperial College London, London SW7 2AZ, UK, and also with the Dept. of Electronic Systems, Aalborg University, Denmark, and with the  Dept.
	of Engineering and Architecture, University of Trieste, Italy.
        E-mail: {\tt\small t.parisini@imperial.ac.uk}}
        \thanks{A. A. Malikopoulos is with the Information and Decision Science Lab, School of Civil $\&$ Environmental Engineering, Cornell University, Ithaca, New York, U.S.A. E-mail: {\tt\small amaliko@cornell.edu}}
        }

\maketitle

\begin{abstract}
This paper investigates the relationship between the team-optimal solution and the Nash equilibrium (NE) to assess the impact of self-interested decisions on team performance. In classical team decision problems, team members typically act cooperatively towards a common objective to achieve a team-optimal solution. However, in practice,  members may behave selfishly by prioritizing their goals, resulting in an NE under a non-cooperative game. To study this misalignment, we develop a parameterized model for team and game problems, where game parameters represent each individual's deviation from the team objective. The study begins by exploring the consistency and deviation between the NE and the team-optimal solution under fixed game parameters. We provide a necessary and sufficient condition for any NE to be a team optimum, along with establishing an upper bound to measure their difference when this consistency fails. { 
We then study how to steer the NE toward the team-optimal solution by adjusting game parameters in an incomplete-information leader--follower setting, where the leader observes only equilibrium responses rather than the exact lower-level game structure. To address this challenge, we develop a two-stage learned-response intervention framework: the leader first learns a behaviorally consistent lower-level model from observed equilibria and then computes the intervention over the learned response map using bilevel hypergradient optimization, followed by convergence analysis and simulation validation.
} 

\end{abstract}


\section{Introduction}
\label{sec:introduction}

Team theory~\cite{marschak1955elements,radner1962team, marschak1972economic} is a mathematical formalism for cooperative decision-making problems in social sciences and engineering. Specifically, a ``team'' consists of several members cooperating to achieve a common objective.  Team theory can be effectively applied in several domains, including intelligent transportation systems~\cite{zhao2019enhanced,Malikopoulos2020}, networked control systems~\cite{hespanha2007survey,zhang2019networked}, and cyber-physical systems~\cite{Malikopoulos2022a}.
The underlying structure of a team decision problem consists of: (i) a team made of a finite number of members; (ii) the decision strategies of each member; (iii) the information available to each member, which may differ across members; (iv) a global objective, which is the same for all members; and (v) the existence, or not, of communication between team members.  
The team outcome results from the joint decisions of all members and is characterized by the team-optimal solution~\cite{radner1962team,zoppoli2020neural,Malikopoulos2021}. This solution refers to a strategy profile under which the overall team performance cannot be improved by changing the strategy of one or more members~\cite{Dave2021a}. 

{Team theory was pioneered by the seminal works~\cite{marschak1955elements,radner1962team, marschak1972economic} on static team problems, and later extended to dynamic settings by~\cite{witsenhausen1971information,witsenhausen1973standard}. 
 In this paper, we focus on static team decision problems. The distinction between static and dynamic teams lies in the underlying information structure that governs team decision-making problems.
In static team problems, the information accessible to each team member remains unaffected by the decisions made by others~\cite{radner1962team,zoppoli2020neural}. Conversely, dynamic team problems involve scenarios in which the information available to at least one member is influenced by the decisions of other members~\cite{Malikopoulos2021}.
}


{Building upon the understanding of 
team problems,
obtaining the team optimum implicitly assumes that all team members are willing to cooperate towards a common objective ~\cite{radner1962team,zoppoli2020neural}.  However,  there may be selfish  members  who prioritize their individual objectives based on perceived interests and preferences
\cite{haurie1985relationship,caporael1989selfishness,altman2007evolutionary}. One reason is that achieving the team optimum often requires members to sacrifice their own benefits, such as tolerating longer transport time in logistics to reduce overall network costs~\cite{shivshankar2014evolutionary}, or consuming additional energy in wireless networks to assist data transmission~\cite{ orda1993competitive}. 
As a result, members may be reluctant to act in accordance with the team goal.
This misalignment
naturally leads to a non-cooperative game formulation~\cite{altman2007evolutionary,haurie1985relationship}, where members make decisions based on the well-known Nash Equilibrium (NE)~\cite{bacsar1998dynamic}. The NE represents a stable strategy profile where no member can benefit from unilaterally changing their strategy, given the strategies of others.}

As each member optimally pursues their individual objectives, the resulting NE may be different  from the team-optimal solution, consequently degrading the overall performance. For instance, in the context of freight transportation~\cite{xu2025deviation}, some shippers may opt for the least expensive routes to maximize their profits. This behavior can lead to congestion and increased costs for other participants, ultimately causing the total system cost to exceed that of the team-optimal solution. Nevertheless,  if the NE aligns with the team optimum, the overall performance of the team remains unchanged. 
Alternatively, if the NE is close to, or can be steered toward the team optimum, the deviation from the desired team outcome remains within an acceptable range and can be further reduced.



{The relationship between  Nash Equilibria (NEs) and team-optimal solutions has been extensively investigated within the framework of team games \cite{bacsar1998dynamic,marden2018game,yuksel2024stochastic}. In a team game setting, all members, considered as players in the game, have identical preferences and aim to minimize or maximize a common objective \cite{marden2018game}. Within the literature on team games, an NE is often referred to as person-by-person optimality in the context of the team problem \cite{yuksel2024stochastic,zoppoli2020neural}. Generally, while any team-optimal solution is guaranteed to be an NE, the converse is not necessarily true. Some studies have delineated conditions under which NE solutions can also be classified as team-optimal in both static and dynamic team games \cite{yuksel2024stochastic}. In the realm of static team games, these conditions are associated with properties such as convexity and Lipschitz continuity. 

{Beyond classical team games, the relationship between NEs and team-optimal solutions in broader game contexts with distinct preferences warrants further exploration. To the best of our knowledge, there has been limited work on this topic.
In this connection, as a practical example, in~\cite{xu2025deviation} we analyzed the consistency and deviation conditions for an NE to be team-optimal within a specific class of flow assignment problems with linear-quadratic objectives.}




In this paper, our primary objective is to examine the effects of self-interested behaviors exhibited by team members on overall team performance. We analyze the consistency and deviation between the NE and team-optimal solution. Additionally, we investigate the design of NE strategies aimed at approaching the team-optimal solution more closely. 

In more specific terms, our contributions are as follows:
\begin{enumerate}
	\item 
We present a comprehensive parameterized model that articulates both cooperative team and non-cooperative game contexts, thereby enabling a structural analysis of the relationship between equilibrium and optimum. The team objective is defined by a set of common parameters, while we represent the individual preferences of each member, in relation to the team preference, through distinct game parameters.

\item We establish the consistency and deviation relationships between team-optimal solutions and NEs under fixed game parameters. {We provide a necessary and sufficient condition indicating that any NE qualifies as a team-optimal solution.}
Additionally, since the consistency condition may not always be satisfied, we establish an upper bound on the deviation of NEs from a team optimum.
\item 
{We address the problem of steering NE strategies toward team optimality by adjusting game parameters in an incomplete-information leader--follower setting, where the leader has access only to observed equilibrium responses rather than an explicit equilibrium-response map. We develop a two-stage learned-response framework: the leader first learns a behaviorally consistent lower-level model from equilibrium data via inverse learning, and then computes the parameter adjustment over the learned response map using an Adam-type hypergradient method.
}

\end{enumerate}

It is worth noting that the scope of this work does not incorporate stochastic aspects. Many practical problems for which our analysis could be useful involve inherent randomness, and exploring such extensions is left for future investigation.

The rest of the paper is organized as follows. After providing some preliminary definitions and concepts in Section~\ref{pre},  Section~\ref{sec: model}     formulates the team and game problems. Section~\ref{sec: Consistency and deviation} then investigates the consistency and deviation relations between the NE and the team-optimal solution.
Subsequently,  Section~\ref{sec:deviation reduction} develops an intervention scheme to guide NE strategies towards  team-optimality, combining an
inverse learning design with a hyper-gradient update.  Section~\ref{sec:diss} provides a comparative discussion, while Section~\ref{sec: simulation} presents numerical experiments to illustrate the theoretical results. Finally, Section~\ref{sec: conclusion} concludes the paper.

\section{Preliminaries}\label{pre}
In this section, we introduce essential notation and provide foundational background material that will be referenced throughout this paper. 

\noindent\textbf{Notation}: 
Let $\mathbb{N}$ denote the set of non-negative integers, $\mathbb{R}$ and $\mathbb{R}_{+}$ denote the sets of real and positive real numbers, respectively,
$ \mathbb{R}^{n} $(or $ \mathbb{R}^{m\times n}, m, n\in\mathbb{N} $) denote the set of $ n $-dimensional (or $ m $-by-$ n $) real column vectors (or real matrices),  $I_n$ denote the $n\times n$ identity matrix,
$ \boldsymbol{1}_{n} $(or $ \boldsymbol{0}_{n} $) denote the $ n $-dimensional column vector with all entries equal to $ {1} $ (or $ {0} $), 
 $ \operatorname{col}\{u_{1}\!,\dots\!,u_{n}\} \!  = \!(\!u^\top_{1},\! \dots,\!u^\top_{n}\!)^\top $, $u_i\in\mathbb{R}^m$, $i=1,\dots,n$,
 $\|\cdot\|$ denote the Euclidean norm of vectors,
$\mathbb{S}^{n}_{\succ 0}$ denote the set of $n\times n$ real symmetric positive definite matrices
 and  $ \nabla f $ denote  the gradient of a differentiable function $f$. For a   differentiable vector-valued mapping $F:\mathbb{R}^n\rightarrow \mathbb{R}^q$, we denote by $\text{J}F(x)\in \mathbb{R}^{q \times n}$  the
 Jacobian of $F$. More generally, 
 for $F:\mathbb{R}^n\times\mathbb{R}^m\rightarrow \mathbb{R}^q$,  we write  $\text{J}_{x}F(x,y)\in \mathbb{R}^{q\times n}$ and $\text{J}_{y}F(x,y)\in \mathbb{R}^{q\times m}$ as the
partial Jacobians of $F$ with respect to its first and second
arguments, respectively. When $q=1$, these reduce to the partial gradients $\nabla_{x}F$ and $\nabla_{y}F$. 

\noindent\textbf{Convex Analysis}:
A set $ K \subseteq \mathbb{R}^{n} $ is convex if $ \omega x_{1}+(1-\omega)x_{2} \in K$ for  any
$ x_{1}, x_{2}\in K $ and $ 0\leq\omega\leq 1   $. 
A mapping $ F : \mathbb{R}^{n} \rightarrow \mathbb{R}^{n} $ is said to be $ \kappa $-strongly monotone on a set $ K $ if there exists a constant $ \kappa > 0 $  such that
$
(F(x)-F(y))^{\top}(x-y) \geq \kappa\|x-y\|^{2}, \quad \forall x, y \in K.
$
The mapping $ F $ is said to be $ \nu $-Lipschitz continuous on set $ K $ if there exists a constant $ \nu > 0 $  such that
$
\|F(x)-F(y)\| \leq \nu\|x-y\|, \quad \forall x, y \in K.
$
For a closed convex set $ K $, the projection operator $ \Pi_{K}: \mathbb{R}^{n} \rightarrow K $ is defined as
$
\Pi_{K}(x) :={\operatorname{argmin}}_{y \in K}\|x-y\|.
$
The projection map is 1-Lipschitz continuous, i.e., $
\|\Pi_{K}(x)-\Pi_{K}(y)\| \leq \|x-y\|,\, \forall x, y \in \mathbb{R}^n.
$ 
Take $  Z\subseteq  \mathbb{R}^{n} $ as a non-empty closed set. 
For $ y \in \mathbb{R}^{n}  $,  denote $ \operatorname{dist}(y,Z) $ as the distance between $ y $ and $ Z $, i.e.,
$  \operatorname{dist}(y,Z) := \inf \limits_{z\in Z}\|y- z\|. 
 $

\noindent\textbf{Non-smooth Analysis}: Any locally Lipschitz mapping
$F:\mathbb{R}^n\rightarrow \mathbb{R}^q$ is almost everywhere differentiable in the sense
of Lebesgue measure by Rademacher’s Theorem \cite{bolte2021conservative}. Its \textit{Clarke Jacobian} is a set-valued mapping $\text{J}^c F : \mathbb{R}^n \rightrightarrows \mathbb{R}^{q \times n}$, defined as
$
\text{J}^c F: = \mathrm{conv}\left\{ \lim_{k \to \infty} \text{J}F(x^k) \;\middle|\; x^k \to x,\; x^k \in \Omega \right\}
$,
where $\Omega \subseteq \mathbb{R}^n$ is the full-measure set of points where $F$ is differentiable. Clarke Jacobians generalize the classical subdifferential from convex analysis to non-smooth, non-convex settings.

The conservative Jacobian generalizes Clarke Jacobians to broader non-smooth settings~\cite{bolte2021conservative}.
Given a locally Lipschitz continuous function $F:\mathbb{R}^n\rightarrow \mathbb{R}^q$, we say that $\mathcal{J}_{F}:\mathbb{R}^n\rightrightarrows \mathbb{R}^{q \times n}$  is a \textit{conservative mapping} or a \textit{conservative Jacobian} for $F$ if $\mathcal{J}_{F}$ has a closed graph, is locally bounded, and is non-empty with
$
\frac{d}{dt}F(x(t))\in \mathcal{J}_{F}(x(t))\dot{x}(t) \; \text{almost everywhere,}
$
where $x:\mathbb{R}\rightrightarrows \mathbb{R}^{n }$ is an absolutely continuous function. When $q=1$, the corresponding vectors are called conservative gradient fields. 
A function $F$ that admits a conservative Jacobian is called \textit{path differentiable}. 
The Clarke and conservative Jacobians satisfy the inclusion $\text{J}^c F(x)\subset \mathcal{J}_{F}(x)$. Unlike Clarke Jacobians, conservative Jacobians allow for operational calculus, e.g., a non-smooth implicit function theorem \cite{bolte2021conservative}.

\section{Cooperative team problems with selfish decision makers}  \label{sec: model}

In this section, we begin with a parameterized framework that captures both the cooperative team and non-cooperative game settings, followed by examples that illustrate its specific forms, applicability, and extensions.

\subsection{Problem formulation}\label{sec: formulation}

Consider a  team problem with $N\in\mathbb{N}$ team members,  indexed by $\mathcal{N}=\{1,\dots,N\}$.  For $i\in\mathcal{N}$,  the $i$-th team member  chooses a team strategy $u_i\in\mathbb{R}^n$ from a {{feasible set}}
$\Xi_i\subseteq \mathbb{R}^n$.  {Throughout this paper, we consider games without coupling constraints among the members, so that the joint feasible set takes the form
$\boldsymbol{\Xi}=\prod_{i=1}^{N}\Xi_i \subseteq \mathbb{R}^{nN}$.}
Let $\boldsymbol{u}=\operatorname{col}\{u_i\}_{i=1}^N\in \boldsymbol{\Xi}$ denote the strategy profile of all members, and let 
$\boldsymbol{u}_{-i}=\operatorname{col}\{u_1,\dots,u_{i-1},u_{i+1},\dots,u_N\}$ denote the strategy profile of all members except member $i$.

All team members independently determine their  strategies and collaborate in minimizing a common cost function $\mathcal{C}(u_1,\dots,u_N):\mathbb{R}^{nN}\rightarrow \mathbb{R}$,  which takes the form of
{\begin{equation}\label{team_cost}
\mathcal{C}(\boldsymbol{u})=\sum\limits_{i=1}^N [f_i(u_i;\alpha)+  g_i(u_i;\beta)^\top\sigma(\boldsymbol{u})+\gamma^\top u_i].
\end{equation}}
   The parameter $\alpha\in\mathbb{R}^{d_\alpha}, d_\alpha\in\mathbb{N},$ in \eqref{team_cost} is associated with each local cost function   $f_i:\mathbb{R}^{d_\alpha}\times\mathbb{R}^{n}\rightarrow \mathbb{R}$. The function $f_i$ represents the cost incurred by member $i$’s own decision in a nonlinear form. The parameter $\beta \in \mathbb{R}^{d_\beta}$, $d_\beta \in \mathbb{N}$, is associated with the cost function $g_i: \mathbb{R}^{d_\beta}\times\mathbb{R}^{n} \rightarrow \mathbb{R}^n$. { The function $\sigma(\boldsymbol{u}):\mathbb{R}^{Nn}\to\mathbb{R}^{n}$ is an aggregative term defined by $$\sigma(\boldsymbol{u})=\sum_{i=1}^N \sigma_i u_i,$$ and $\sigma_i\in\mathbb{R}$ corresponds to the weight associated with member $i$ with $\boldsymbol{\sigma}=(\sigma_1,\dots, \sigma_N)\in \mathbb{R}^N$. The function $g_i$, through its coupling with $\sigma(\boldsymbol{u})$, characterizes how the other members influence member $i$'s cost, capturing effects such as congestion, interference, or shared resource usage.} The linear term $\gamma^\top u_i$ reflects team-level incentives or penalties, where $\gamma \in \mathbb{R}^{d_\gamma}$, $d_\gamma =n$ denote the corresponding cost coefficients. {
   This formulation emphasizes the interdependencies between individual and collective decisions within the model.
   }

 This unified formulation finds applications in various scenarios, such as traffic routing \cite{burton1992instance,allen2022using} and wireless communication \cite{sung1999power}. For instance, in a shortest path traffic problem \cite{burton1992instance,allen2022using}, team members, viewed as vehicles operated by the same carrier, allocate flow across network arcs to minimize the total travel time. 
The travel time includes a linear term accounting for the free-flow travel time and a quadratic term capturing the additional delay caused by congestion from the aggregated flow \cite{allen2022using}.  We will present two specific examples in Section \ref{example}.

The team problem under consideration follows a static information structure. That is,  the information available for each member’s decision-making, such as the cost function and feasible set constraints,
is not influenced by the decisions of others~\cite{Malikopoulos2021,zoppoli2020neural}. 
The team decision-making problem is formulated as follows:
\begin{equation}\label{team_function}
\begin{aligned}
&\min_{u_1,\dots,u_N}\; \mathcal{C} (u_1,\dots,u_N) \\
& \mathrm{s.\,t.}
\;  
u_i\in \Xi_i,\;
i\in\mathcal{N}.
\end{aligned}
\end{equation}
{
The team problem \eqref{team_function} is a one‑shot (single‑stage) problem,
where each member takes its entire decision strategy $u_i$ at once. In fact, the formulated model admits a natural extension to the multi‑stage setting.
For instance, each member chooses a control sequence  \(\{u_{i,t}\}_{t=0}^{T-1}\) over \(T\) time steps, and the team objective is a finite‑horizon  linear–quadratic regulator (LQR) cost. We will also give details in Section \ref{example}.}

The corresponding outcome  of \eqref{team_function} is described by a team-optimal solution \cite{zoppoli2020neural,Malikopoulos2021,yuksel2024stochastic}.

\begin{mydef}[Team-optimal solution]
A strategy profile $\boldsymbol{u}^*=\operatorname{col}\{u_1^*, u_2^*, \dots, u_N^*\}\in \boldsymbol{\Xi}$ is a team-optimal solution if 
\begin{align*}
    \mathcal{C}(\boldsymbol{u}^*)\leq  \mathcal{C}(\boldsymbol{u}), \; \forall \boldsymbol{u} \in \boldsymbol{\Xi}.
\end{align*}
\end{mydef}
\begin{remark}
The term ``strategy'' in this paper simply refers to a decision variable, consistent with its usage in game theory. Unlike in classical team theory \cite{zoppoli2020neural}, we do not explicitly formalize strategies as functions of information,  as the role of information is not the focus of our study.
\end{remark}

The team-optimal solution $\boldsymbol{u}^*$ can be obtained either in a centralized or a decentralized way. Following \cite{Malikopoulos2021}, the decision-making process can be viewed from the perspective of a ``team manager'' who has full prior knowledge of the information of all members and seeks to calculate $\boldsymbol{u}^* $ in a centralized way.
Alternatively,
in a decentralized setting, each team member makes decisions based on their own information, and estimates others’ strategies through communication, ultimately achieving $\boldsymbol{u}^* $ collaboratively.  In both cases, the resulting strategy of each team member is the same as the one derived by the manager \cite{Malikopoulos2021}. 

 However, such alignment implicitly assumes that all team members are committed to pursuing a common objective $\mathcal{C}(\boldsymbol{u})$ \cite{radner1962team,Malikopoulos2021}. In practice, members may instead prioritize minimizing their own objectives due to self-interest \cite{altman2007evolutionary,haurie1985relationship,xu2025deviation}. This leads to  a non-cooperative game formulation, where each member has an individual cost function $\mathcal{C}_i(u_{i}, \boldsymbol{u}_{-i}):\mathbb{R}^{nN}\rightarrow \mathbb{R} $, given by
{\begin{equation}\label{para}
  \mathcal{C}_i(u_i,\boldsymbol{u}_{-i})= f_i(u_i;\alpha_i)+  g_i(u_i;\beta_i)^\top\sigma(\boldsymbol{u})+\gamma_i^\top u_i,
\end{equation}}
where the parameters $\alpha_i\in \mathbb{R}^{d_\alpha}$, $\beta_i\in \mathbb{R}^{d_\beta}$ and $\gamma_i\in \mathbb{R}^n$ reflect individual preferences,
in contrast to the team-level parameters $\alpha$, $\beta$ and $\gamma$ in \eqref{team_cost}. The cost function $\mathcal{C}_i$ can be viewed as a personalized decomposition of the team cost $\mathcal{C}$, in which each member focuses solely on its local cost $f_i$,  interaction term $g_i$, and linear term $\gamma_i^\top u_i$. When $\alpha_i=\alpha$, $\beta_i=\beta$, and $\gamma_i=\gamma$ for all $i\in\mathcal{N}$, we have $\sum_{i=1}^N\mathcal{C}_i(\boldsymbol{u})=\mathcal{C}(\boldsymbol{u})$.


Given $\boldsymbol{u}_{-i}$, team member $i$ aims to solve
\begin{equation}\label{f1}
	\min_{u_{i} \in \Xi_{i}} \;\mathcal{C}_i\left(u_{i}, \boldsymbol{u}_{-i}\right). \quad  
\end{equation}
The outcome of  game  \eqref{f1} is characterized by the NE \cite{nash1951non}.

\begin{mydef}[NE]\label{d2}
	A profile $ \boldsymbol{u}^{\Diamond}=\operatorname{col}\{u_{1}^{\Diamond}, \dots, u_{N}^{\Diamond} \} \in \boldsymbol{\Xi}
    $ is said to be an  NE 
    if 
	\begin{equation*}\label{ne}
		\mathcal{C}_i\left(u_{i}^{\Diamond}, \boldsymbol{u}_{-i}^{\Diamond}\right) \leq \mathcal{C}_i\left(u_{i}, \boldsymbol{u}_{-i}^{\Diamond}\right),\;\forall i \in \mathcal{N}, \forall u_{i}\in \Xi_{i}. 
	\end{equation*}
\end{mydef}
{The NE describes an outcome where no team member has an incentive to unilaterally deviate from their chosen strategy, provided that the strategies of all others remain unchanged.}

\begin{remark}
When all team members have identical preferences, i.e., $\mathcal{C}_1=\dots=\mathcal{C}_N$, the game becomes a team game \cite{marden2018game}. In this setting, an NE reflects person‑by‑person optimality within the team context, which is a weaker solution concept than a team-optimal solution \cite{zoppoli2020neural,yuksel2024stochastic}.
Although any team-optimal solution qualifies as an NE, the converse does not necessarily hold. We will further elaborate on this specific case in the subsequent section.
\end{remark}

Given that individual objectives diverge from the team goal, the resulting NE may fail to be team-optimal and thus compromise overall team performance. It is of interest to explore when NE strategies coincide with the team-optimal solution.
 We seek to address the following questions:
\begin{itemize}
      \item Under what conditions are the team-optimal solution and the NE consistent?
    \item Does there exist a deviation bound between the team-optimal solution and the NE if consistency fails?
    \item 
    Can this deviation be mitigated by adjusting game parameters to influence members' preferences?
\end{itemize}

Now, we provide two examples that illustrate the specific forms and applicability of \eqref{team_cost} and \eqref{para}, followed by a third example for the model extension.
}

\subsection{Examples}\label{example}

\noindent \textbf{Example 1}  Consider a traffic network with $N$ vehicles belonging to the same carrier, where each vehicle $i\in\mathcal{N}$ allocates flow  $u_i \in \mathbb{R}^n$ across the network arcs \cite{burton1992instance,allen2022using}. 
The team objective is to minimize the total travel time of all vehicles, which depends quadratically on the aggregate flow across the network and includes free-flow travel costs.
\begin{align}\label{team_traffic}
\mathcal{C}(\boldsymbol{u})= \sum\limits_{i=1}^N u_i^\top Q(\alpha)  u_i + u_i^\top Q(\beta) \sigma(\boldsymbol{u})+\gamma^\top u_i.
\end{align}
Here, $Q(\alpha)\in \mathbb{R}^{n\times n}$ and $Q(\beta)\in \mathbb{R}^{n\times n}$ are positive diagonal matrices. In particular,  $Q(\alpha)=\sum_{l=1}^{d_\alpha}Q_{l}[\alpha]_l$, where $[\alpha]_l$ denotes the $l$-th component of $\alpha \in\mathbb{R}^{d_\alpha}$, and  $Q(\beta)$ is defined similarly. 
The matrix $Q(\alpha)$  captures the coefficients for the individual cost associated with each vehicle’s own flow, while $Q(\beta)$ represents the interaction cost coefficients between each vehicle and the aggregate flow of others across the network. 
The vector 
$\gamma\in\mathbb{R}^n$ represents the base cost for traveling each arc. In contrast, each vehicle is only concerned with the cost induced by its own flow and the congestion caused by others on the same routes.
Its individual cost function is 
\begin{align}\label{game_traffic}
\mathcal{C}_i(u_i,\boldsymbol{u}_{-i})\!=\! u_i^\top Q(\alpha_i)  u_i \!+\! u_i^\top Q(\beta_i) \sigma(\boldsymbol{u})\!+\!\gamma_i^\top u_i,
\end{align}
where $Q(\alpha_i)  $, $Q(\beta_i)$, and $\gamma_i$ represent vehicle
$i$'s subjective preferences of travel costs.

{\noindent \textbf{Example 2} Consider a multi-channel wireless uplink power-allocation
problem with \(N\) users transmitting over \(n\) independent subchannels.
Each user \(i\in\mathcal N\) chooses a power-allocation vector
\(u_i=\operatorname{col}\{u_{i\ell}\}_{\ell=1}^n\in\mathbb R^n\),
where \(u_{i\ell}\) denotes the transmit power of user \(i\) on subchannel
\(\ell\).
For each subchannel \(\ell\), define the aggregate interference/load term
$[\sigma(\boldsymbol u)]_\ell=\sum_{k=1}^N \sigma_k u_{k\ell}$,
$\sigma_k>0$.
The team objective is given by
\[
\mathcal{C}(\boldsymbol{u})
\!=\!\!
\sum_{i=1}^N
\sum_{\ell=1}^n
-\alpha_\ell \log\bigl(1+h_{i\ell}u_{i\ell}\bigr)
\!+\!
\beta_\ell \sigma_i u_{i\ell}[\sigma(\boldsymbol u)]_\ell
\!+\!
\gamma_\ell u_{i\ell},
\]
where \(h_{i\ell}>0\) is the channel gain of user \(i\) on subchannel
\(\ell\), \(\alpha_\ell>0\) weights the rate benefit, \(\beta_\ell>0\)
weights the interference-pricing cost, and \(\gamma_\ell>0\) is the
linear energy cost coefficient. In contrast, user \(i\) acts selfishly according to
\[
\mathcal{C}_i(u_i,\boldsymbol{u}_{-i})
\!=\!\!
\sum_{\ell=1}^n\!
-\alpha_{i\ell} \log\bigl(1\!+\!h_{i\ell}u_{i\ell}\bigr)
\!+\!
\beta_{i\ell}\sigma_i u_{i\ell}[\sigma(\boldsymbol u)]_\ell
\!+\!
\gamma_{i\ell}u_{i\ell}
,
\]
where \(\alpha_{i\ell}\), \(\beta_{i\ell}\), and \(\gamma_{i\ell}\)
represent user \(i\)'s subjective rate-benefit, interference-pricing, and
energy-cost parameters on subchannel \(\ell\).

}

\noindent \textbf{Example 3}  Consider a finite-horizon linear--quadratic regulator (LQR) problem in which each member \(i\) controls a sequence of inputs \(\{u_{i,t}\in\mathbb{R}\}\) over \(T\) time steps.  The joint state evolves according to 
$x_{t+1}=Ax_t+B\sum_{i=1}^N \sigma_i u_{i,t}$, $t=0,\dots,T-1$. The team seeks to minimize:
$\sum_{t=0}^{T-1}(x_t^\top Q(\tilde{\alpha})x_t+\sum_{i=1}^N u_{i,t}^\top R_i(\tilde{\beta})u_{i,t})+x_T^\top Q_f(\tilde{\alpha})x_T$,
where the per-step cost weights \(Q(\tilde{\alpha})\), \(Q_f(\tilde{\alpha})\), and \(R_i(\tilde{\beta})\) are parameterized by \(\tilde{\alpha}\) and \(\tilde{\beta}\). {By stacking the sequence of controls \(\{u_{i,t}\}_{t=0}^{T-1}\) into a single vector \(u_i\in\mathbb{R}^T\), one can
rewrite the multi-stage problem as the one-shot optimization
\[
\min_{u_1,\dots,u_N}
\sum_{i=1}^N(
\tfrac12\,u_i^\top H_i(\tilde{\alpha},\tilde{\beta})u_i
+
u_i^\top G_i(\tilde{\alpha},\tilde{\beta})\,\sigma(\boldsymbol u)
+
\gamma^\top u_i
),
\]
where \(H_i(\tilde{\alpha},\tilde{\beta})\) collects the quadratic self-cost terms of member \(i\), and \(G_i(\tilde{\alpha},\tilde{\beta})\) captures the coupling induced by the shared dynamics and state penalties.}

\section{Consistency and deviation between NE and  team-optimal solution}\label{sec: Consistency and deviation}


In this section, we address the first two questions posed in Section~\ref{sec: formulation}. Specifically, we study when NEs and team-optimal solutions coincide, and how far they may differ when consistency fails.

\subsection{Consistency relation}
We first impose the following assumptions on Team Problem~\eqref{team_function} and Game~\eqref{f1} to study the consistency relation.
\begin{assumption}\label{assum2} 
	$\ $
	\begin{enumerate}[(1)]
    \item 
    For each $i\in\mathcal{N}$, the set ${\Xi}_i$ is a nonempty closed convex polyhedron  of the form
    $
{\Xi}_i=\{u_i\in\mathbb{R}^n|D_i u_i\leq b_i, H_i u_i=m_i\},
    $
    where $D_i\in\mathbb{R}^{q_i\times n}, q_i\in\mathbb{N},$ and $H_i\in\mathbb{R}^{r_i\times n}, r_i\in\mathbb{N},$ are full row rank, and $b_i\in\mathbb{R}^{q_i}$, $m_i\in\mathbb{R}^{r_i}$. The set ${\Xi}_i$  satisfies Slater’s constraint qualification.
 \item    The function $\mathcal{C}(\boldsymbol{u})$ is convex
and continuously differentiable in $\boldsymbol{u}$. 
For each  $i\in\mathcal{N}$, the function $\mathcal{C}_i(u_i,\boldsymbol{u}_{-i})$ is convex in  $u_i$
and continuously differentiable in $\boldsymbol{u}$. 
\item For each $i\in\mathcal{N}$, the set ${\Xi}_i$ is compact.
\end{enumerate}
\end{assumption}

Assumption \ref{assum2} guarantees the existence of NEs and team-optimal solutions, without requiring uniqueness. 
{ The following results verify the existence of a team-optimal solution and an NE \cite[Proposition 1.4.2, Proposition 2.2.9]{facchinei2003finite}. }
	\begin{lemma}\label{l1}
			Given Assumption \ref{assum2},  there exists a team-optimal solution  $\boldsymbol{u}^{*}$ in  team decision problem \eqref{team_function}.
		\end{lemma}
\begin{lemma}\label{l2}
			Given Assumption \ref{assum2},  there exists an NE  $\boldsymbol{u}^{\Diamond}$ in game \eqref{f1}.
		\end{lemma}

Based on Lemmas \ref{l1} and \ref{l2}, the following result establishes the equivalence between NEs and team-optimal solutions in two cases. The first  case corresponds to the classic team game, in which all members have identical cost functions. {The second case is related to potential games~\cite{monderer1996potential}, where the members' cost functions satisfy a potential-type consistency condition:
\begin{equation}\label{potential}
   \mathcal{C}(u_{i}^{\prime}, \boldsymbol{u}_{-i})-\mathcal{C}\left(\boldsymbol{u}\right)
   =
   \mathcal{C}_i(u_{i}^{\prime}, \boldsymbol{u}_{-i})-\mathcal{C}_i\left(\boldsymbol{u}\right), 
\end{equation}        
for any $\boldsymbol{u}\in \boldsymbol{\Xi}$, $i\in\mathcal{N}$, and unilateral deviation $u_i^{\prime}\in \Xi_i$.
Condition~\eqref{potential}  ensures that the unilateral variation of each member's cost is exactly captured by the variation of the prescribed team cost. 
}


\begin{proposition}\label{p1} 
Suppose that either (a) all team members have the same cost function $\mathcal{C}$  as in \eqref{team_cost}, or (b) their cost functions satisfy the condition in \eqref{potential}.  Then: (i) Any  team-optimal solution $\boldsymbol{u}^*$ is an  NE $\boldsymbol{u}^\Diamond$; (ii) Under Assumption \ref{assum2}, any NE $\boldsymbol{u}^\Diamond$ is also  a team-optimal solution $\boldsymbol{u}^*$.
\end{proposition}

\noindent\textbf{Proof.} Conclusion (i) follows directly for cases (a) and (b) from the definition of the team-optimal solution.

We now proceed to the proof of conclusion~(ii), beginning with case~(a). 
Let \( \boldsymbol{u}^\Diamond \) be an NE. Since \( \Xi_i \) is compact and convex, and the team cost function \( \mathcal{C} \) is convex and continuously differentiable in \( u_i \), the first-order optimality condition (by convexity and the minimum principle) implies that
\[
(u_i - u_i^\Diamond)^\top \nabla_{u_i} \mathcal{C}(\boldsymbol{u}^\Diamond) \geq 0, \quad \forall u_i \in \Xi_i, \forall i \in \mathcal{N}.
\]
By stacking these inequalities across all \( i \in \mathcal{N} \), we obtain
$(\boldsymbol{u} - \boldsymbol{u}^\Diamond)^\top \nabla \mathcal{C}(\boldsymbol{u}^\Diamond) \geq 0$,  $\forall \boldsymbol{u} \in \boldsymbol{\Xi}$,
which is the first-order condition for optimality of the team objective \( \mathcal{C} \) over the feasible set \( \boldsymbol{\Xi} \). Hence, \( \boldsymbol{u}^\Diamond \) is also a team-optimal solution, satisfying
$\mathcal{C}(\boldsymbol{u}^\Diamond) \leq \mathcal{C}(\boldsymbol{u}), \;\forall \boldsymbol{u} \in \boldsymbol{\Xi}.$
This concludes the first part of the proof.

As for case~(b), since the functions \( \mathcal{C}_i \) are continuously differentiable, condition~\eqref{potential} is equivalent to
\[
\nabla_{u_i} \mathcal{C}_i(u_i, \boldsymbol{u}_{-i}) = \nabla_{u_i} \mathcal{C}(u_i, \boldsymbol{u}_{-i}) \quad \text{for } i \in \mathcal{N}.
\]
Under this condition, the gradient structure of the game aligns with that of the team cost function. Hence, the first-order conditions satisfied by the NE for the individual objectives \( \mathcal{C}_i \) are identical to those for the team objective \( \mathcal{C} \). The proof for case~(b) then follows analogously to case~(a).
\hfill $\square$



In general, this equivalence does not hold. We next characterize when an NE is team-optimal beyond the two special cases above. Consider the partial gradient of $\mathcal{C}$ in \eqref{team_cost}  in $u_i$, given by
\begin{align*}
\nabla_{u_i}\mathcal{C}(u_i,\boldsymbol{u}_{-i})=&\nabla_{u_i}f_i(u_i;\alpha)+ \nabla_{u_i}g_i(u_i;\beta)^\top\sigma(\boldsymbol{u})\\
&\,+\sigma_i\!\sum\nolimits_{j=1}^N g_j(u_j;\beta)+\gamma, 
\end{align*}
and the  partial gradient of $\mathcal{C}_i$ in \eqref{para} with respect to $u_i$, 
\begin{align}\label{gradientne}
\nabla_{u_i}\mathcal{C}_i(u_i,\boldsymbol{u}_{-i})=&\nabla_{u_i}f_i(u_i;\alpha_i)+ \nabla_{u_i}g_i(u_i;\beta_i)^\top\sigma(\boldsymbol{u})\notag\\
&+\sigma_ig_i(u_i;\beta_i)+\gamma_i.
\end{align}
Let {$\operatorname{int}(\Xi_{i})$ denote the interior of the set $\Xi_{i}$} and $\delta(\cdot)$ denote  a neighbourhood of a  point in $\mathbb{R}^n$.  
The following result establishes the condition for any NE to be team-optimal when the strategy sets are one-dimensional{\footnote{Here, we exclude a  degenerate case in which the individual marginal cost \(\nabla_{u_i}\mathcal{C}_i(u_i,\boldsymbol{u}_{-i})\) is identically zero over a nontrivial feasible interval.
}}.


\begin{mythm}\label{t111}
{Under Assumption \ref{assum2},  suppose that the decision variables are one-dimensional so that $\Xi_i=[u_i^{\min},u_i^{\max}]$,  $i\in\mathcal N$. 
Then every  NE $\boldsymbol{u}^\Diamond\in \boldsymbol{\Xi} \subseteq \mathbb{R}^N$  of \eqref{f1}  is a team-optimal solution  of \eqref{team_function}  if and only if, 
for $i\in\mathcal{N}$, there exists  a neighborhood $ \delta(u_i^{\Diamond})$   such that 
    \begin{enumerate}[(i)]
       \item   $ \nabla_{u_i}\mathcal{C}(u_i^{\Diamond},\boldsymbol{u}_{-i}^{\Diamond})=0 $, or  \item $\nabla_{u_i}\mathcal{C}(u_i,\boldsymbol{u}_{-i}^{\Diamond})\cdot \nabla_{u_i}\mathcal{C}_i(u_i,\boldsymbol{u}_{-i}^{\Diamond}) >0 $ for $  u_i \in \delta(u_i^{\Diamond})\cap \operatorname{int}(\Xi_{i})$.
    \end{enumerate}}	
\end{mythm}

\noindent\textbf{Proof.} 
Before presenting the detailed proof, we outline the main ideas through the following cases:
\begin{enumerate}[i)]
    \item The NE strategies lie in the interiors of their strategy sets, and the corresponding partial derivatives are zero.

    \item The NE strategies lie on the boundaries of their strategy sets, and the corresponding partial derivatives are zero.
    
    \item The NE strategy lies on the boundaries of the strategy sets, and the corresponding partial derivatives are nonzero.
\end{enumerate}
{For each case, we show that any NE $\boldsymbol{u}^\Diamond$ of game~\eqref{f1} is also an NE of the team  problem~\eqref{team_function}. Then, by Proposition \ref{p1}(a) (ii), we conclude that such an NE of \eqref{team_function} is  a team-optimal solution $\boldsymbol{u}^*$, thereby completing the proof.}

\textit{Sufficiency}: For each \( i \in \mathcal{N} \), we consider two cases depending on the value of \( \nabla_{u_i} \mathcal{C}(u_i^\Diamond, \boldsymbol{u}_{-i}^\Diamond) \).

- {Suppose \( \nabla_{u_i} \mathcal{C}(u_i^\Diamond, \boldsymbol{u}_{-i}^\Diamond) = 0 \). By the convexity of \( \mathcal{C} \) in \( \boldsymbol{u} \),
$u_i^\Diamond \in \operatorname{argmin}_{u_i \in \Xi_i} \mathcal{C}(u_i, \boldsymbol{u}_{-i}^\Diamond)$, 
which implies that $u_i^\Diamond$ is an NE strategy of the team game associated with cost $\mathcal C$ \eqref{team_function}. Then, by Proposition~\ref{p1}(a) (ii), \( u_i^\Diamond \) is also a team-optimal strategy.}

- {Suppose \( \nabla_{u_i} \mathcal{C}(u_i^\Diamond, \boldsymbol{u}_{-i}^\Diamond) \neq 0 \). By the coincidence condition (ii), we have
$\nabla_{u_i} \mathcal{C}(u_i, \boldsymbol{u}_{-i}^\Diamond) \cdot \nabla_{u_i} \mathcal{C}_i(u_i, \boldsymbol{u}_{-i}^\Diamond) > 0$,  $\forall u_i \in \delta(u_i^\Diamond) \cap \operatorname{int}(\Xi_i)$.
If \( u_i^\Diamond \in \operatorname{int}(\Xi_i) \), then the definition of  NE requires \( \nabla_{u_i} \mathcal{C}_i(u_i^\Diamond, \boldsymbol{u}_{-i}^\Diamond) = 0 \), which contradicts the strict positivity above. Therefore, \( u_i^\Diamond \notin \operatorname{int}(\Xi_i) \), and must lie on the boundary. 
Assume \( u_i^\Diamond = u_i^{\min} \). Then, there exists a neighborhood \( \delta'(u_i^{\min}) \) such that
$\nabla_{u_i} \mathcal{C}_i(u_i, \boldsymbol{u}_{-i}^\Diamond) \geq 0$,  $\forall u_i \in \delta'(u_i^{\min}) \cap \Xi_i$,
and in particular,
$\nabla_{u_i} \mathcal{C}_i(u_i, \boldsymbol{u}_{-i}^\Diamond) > 0$, $\forall u_i \in \delta'(u_i^{\min}) \cap \operatorname{int}(\Xi_i)$.
Define \( \delta''(u_i^{\min}) := \delta(u_i^{\min}) \cap \delta'(u_i^{\min}) \). Then,
$\nabla_{u_i} \mathcal{C}(u_i, \boldsymbol{u}_{-i}^\Diamond) > 0$, $ \forall u_i \in \delta''(u_i^{\min}) \cap \operatorname{int}(\Xi_i)$.
By the continuity and monotonicity of \( \nabla_{u_i} \mathcal{C}(u_i, \boldsymbol{u}_{-i}^\Diamond) \), we conclude that
$\nabla_{u_i} \mathcal{C}(u_i, \boldsymbol{u}_{-i}^\Diamond) > 0, \quad \forall u_i \in \operatorname{int}(\Xi_i)$. 
Therefore,
$u_i^{\min} \in \operatorname{argmin}_{u_i \in \Xi_i} \mathcal{C}(u_i, \boldsymbol{u}_{-i}^\Diamond)$,
i.e., \( u_i^\Diamond = u_i^{\min} \) is an NE strategy of team \eqref{team_function}, and moreover it is team-optimal by Proposition~1 (a) (ii). The same reasoning applies if \( u_i^\Diamond = u_i^{\max} \).}

Hence, if, for each $i$, at least one of (i) or (ii) holds, then any NE \( \boldsymbol{u}^\Diamond \) of game~(4) is team-optimal of \eqref{team_function}.

\textit{Necessity}: Suppose \( \boldsymbol{u}^\Diamond \) is both an NE of~\eqref{f1} and a team-optimal solution of~\eqref{team_function}. Then for each \( i \in \mathcal{N} \):

- If \( u_i^\Diamond \!\in \!\operatorname{int}(\Xi_i) \), \( \nabla_{u_i} \!\mathcal{C}(u_i^\Diamond, \!\boldsymbol{u}_{-i}^\Diamond) \!= \!0 \), which satisfies the condition (i).

- If \( u_i^\Diamond \!=\! u_i^{\min} \), and \( \nabla_{u_i} \!\mathcal{C}(u_i^{\min}, \!\boldsymbol{u}_{-i}^\Diamond) \!=\! 0 \), the condition (i) is satisfied.

Otherwise, we have \( \nabla_{u_i} \mathcal{C}(u_i^{\min}, \boldsymbol{u}_{-i}^\Diamond) > 0 \). Then there exists a neighborhood \( \delta'(u_i^{\min}) \) such that
$\nabla_{u_i} \mathcal{C}(u_i, \boldsymbol{u}_{-i}^\Diamond) > 0, \quad \forall u_i \in \delta'(u_i^{\min}) \cap \operatorname{int}(\Xi_i)$. 
Moreover, since \( u_i^\Diamond = u_i^{\min} \) is an NE strategy, it satisfies
$\mathcal{C}_i(u_i^{\min}, \boldsymbol{u}_{-i}^\Diamond) < \mathcal{C}_i(u_i, \boldsymbol{u}_{-i}^\Diamond), \quad \forall u_i \in \delta''(u_i^{\min}) \cap \operatorname{int}(\Xi_i)$, 
for some neighborhood \( \delta''(u_i^{\min}) \), which implies
$\nabla_{u_i} \mathcal{C}_i(u_i, \boldsymbol{u}_{-i}^\Diamond) > 0, \quad \forall u_i \in \delta''(u_i^{\min}) \cap \operatorname{int}(\Xi_i)$. 
Let \( \delta(u_i^{\min}) := \delta'(u_i^{\min}) \cap \delta''(u_i^{\min}) \). Then
$\nabla_{u_i} \mathcal{C}(u_i, \boldsymbol{u}_{-i}^\Diamond) \cdot \nabla_{u_i} \mathcal{C}_i(\cdot, \boldsymbol{u}_{-i}^\Diamond) > 0, \quad \forall u_i \in \delta(u_i^{\min}) \cap \operatorname{int}(\Xi_i)$.
A similar argument applies for the case \( u_i^\Diamond = u_i^{\max} \). This concludes the proof.
 \hfill $\square$

{ Theorem~\ref{t111} addresses the question of when selfish behavior remains consistent with the original team performance once cooperation is no longer followed. This question has not been sufficiently studied in the existing literature, which typically treats team optimality and equilibrium behavior separately. Theorem~\ref{t111} fills this gap by establishing a necessary and sufficient condition based on local alignment between the marginal individual incentive and the marginal team incentive around the equilibrium. Thus, it explains when selfish behavior is harmless to the team objective and when inefficiency arises, depending on whether the private and team marginal incentives align or fail to do so.
}


When the NE is unique, the following conclusion follows directly from Theorem \ref{t111}.

\begin{corollary}\label{c1}
Under Assumption~\ref{assum2}, suppose the NE $\boldsymbol{u}^\Diamond \in \boldsymbol{\Xi} \subseteq \mathbb{R}^N$ is unique. Then, $\boldsymbol{u}^\Diamond$ is a team-optimal solution to~\eqref{team_function} if and only if, for each $i \in \mathcal{N}$, there exists a neighborhood $\delta(u_i^\Diamond)$ such that
$\nabla_{u_i} \mathcal{C}(u_i, \boldsymbol{u}_{-i}^\Diamond) \cdot \nabla_{u_i} \mathcal{C}_i(u_i, \boldsymbol{u}_{-i}^\Diamond) \geq 0 \quad \text{for } u_i \in \delta(u_i^\Diamond) \cap \Xi_i$.
\end{corollary}

{
\begin{remark}
Theorem~\ref{t111} naturally extends to certain multi-dimensional settings in which the decision variables are coordinatewise decoupled. In particular, when \(u_i\in\mathbb{R}^n\), the feasible set is box-constrained, and both the team and individual costs are separable across the coordinates of \(u_i\), the same local marginal-alignment argument can be applied coordinatewise. Such a setting arises, for instance, in the multi-channel power control problem in Example~2, where each coordinate \(u_{i\ell}\), \(\ell=1,\dots,n\), represents the transmit power on one independent subchannel.
In more general multi-dimensional polyhedral settings, the coordinates may be coupled through the objective or the constraints. In particular, when the profile lies on the boundary of a polyhedral feasible set, the partial derivatives may have mixed signs across coordinates, making a componentwise analysis cumbersome. In such cases, the coincidence relation is more naturally analyzed through variational-inequality conditions.
\end{remark}
}

\vspace{-0.4cm}

\subsection{Deviation bound}

In general, the team-optimal solution and the NE do not
coincide. To assess the reasonableness of individual decisions, we focus on the case of a unique team‐optimal solution and characterize the deviation of the NE set from that solution.

Specifically, we employ the Hausdorff metric to quantify the deviation between sets~\cite{rockafellar2005variational,birsan2005one}, a widely used measure of set-to-set distance that captures the greatest distance a point in one set must travel to reach the other and thereby characterizes the worst-case mismatch. 
For two non-empty {closed} sets \( X, Z \subseteq \mathbb{R}^n \), the Hausdorff metric is defined by
\[
d_H(X, Z) := \max\left\{ \sup_{x \in X} \inf_{z \in Z} \|x - z\|,\ \sup_{z \in Z} \inf_{x \in X} \|x - z\| \right\}.
\]
Accordingly, let $\Upsilon_{TO}$ be the set of  team-optimal strategy profile $\boldsymbol{u}^{*}$ and $\Upsilon_{NE}$ be the set of  NE strategy profile $\boldsymbol{u}^{\Diamond}$. 
Under the convexity and continuity conditions specified in Assumption~\ref{assum2}, both the team-optimal solution and the NE correspond to first-order critical points of their respective optimization problems. This equivalence enables the analysis of their deviation using gradient-based information. Define  
$G(\boldsymbol{u}):=\nabla_{\boldsymbol{u}}\mathcal{C}(\boldsymbol{u})=\operatorname{col}\{\nabla_{{u_1}}\mathcal{C}(\cdot,u_{-1}),\dots, \nabla_{{u_N}}\mathcal{C}(\cdot,u_{-N})\}$
as the gradient map of cost function $\mathcal{C}$  and  
$F(\boldsymbol{u}):=\operatorname{col}\{\nabla_{{u_1}}\mathcal{C}_1(\cdot,u_{-1}),\dots,\nabla_{{u_N}}\mathcal{C}_N(\cdot,u_{-N})\}$
as the
pseudo-gradient of $\mathcal{C}_i$ for $i\in \mathcal{N}$. 
{ Recalling the expressions of $F$ and $G$, we introduce the following Lipschitz continuity assumption on the parameter dependence of \(\nabla_{u_i}f_i\), \(\nabla_{u_i}g_i\), and \(g_i\).
\begin{assumption}\label{assum_param_lip}
For each \(i\in\mathcal N\), there exist constants \(L_i^{\nabla f}>0\), \(L_i^{\nabla g}>0\), and \(L_i^g>0\) such that, for all \(u_i\in \Xi_i\), \(\alpha,\alpha_i\), and \(\beta,\beta_i\),
$\|\nabla_{u_i}f_i(u_i;\alpha_i)-\nabla_{u_i}f_i(u_i;\alpha)\|
\le
L_i^{\nabla f}\|\alpha_i-\alpha\|$,
$\|\nabla_{u_i}g_i(u_i;\beta_i)-\nabla_{u_i}g_i(u_i;\beta)\|
\le
L_i^{\nabla g}\|\beta_i-\beta\|$,
and
$\|g_i(u_i;\beta_i)-g_i(u_i;\beta)\|
\le
L_i^g\|\beta_i-\beta\|$.
\end{assumption}

On this basis, we obtain an a priori bound on the deviation of NE strategies.
}



{\begin{mythm}\label{t13}
Suppose that Assumptions~\ref{assum2} and \ref{assum_param_lip} hold and that there exists a constant \(\kappa_1>0\) such that the mapping \(G(\boldsymbol{u})\) is \(\kappa_1\)-strongly monotone.
Then, 
\[
d_H(\Upsilon_{TO}, \Upsilon_{NE}) \le \frac{ \xi}{\kappa_1},
\vspace{-0.2cm}
\]
where 
\begin{align*} \xi :=& ( \sum_{i=1}^N ( L_i^{\nabla f}\|\alpha_i-\alpha\| + (\bar\sigma\,L_i^{\nabla g}+|\sigma_i|L_i^g)\|\beta_i-\beta\|\\ &\qquad\qquad\qquad + |\sigma_i|\sum_{j\neq i}\bar g_j + \|\gamma_i-\gamma\| )^2 )^{1/2}, \end{align*}
and \(\bar\sigma>0\) and \(\bar g_j>0\) are constants.

\end{mythm}

\noindent\textbf{Proof.}  
Take any \(\boldsymbol{u}^{\Diamond}\in \Upsilon_{NE}\) and \(\boldsymbol{u}^{*}\in \Upsilon_{TO}\). Under Assumption 1, since \(\boldsymbol{u}^{\Diamond}\) is an NE, it satisfies \(F(\boldsymbol{u}^{\Diamond})^\top(\boldsymbol{u}-\boldsymbol{u}^{\Diamond})\ge 0\) for all \(\boldsymbol{u}\in \boldsymbol{\Xi}\). By choosing \(\boldsymbol{u}=\boldsymbol{u}^{*}\), we obtain \(F(\boldsymbol{u}^{\Diamond})^\top(\boldsymbol{u}^{*}-\boldsymbol{u}^{\Diamond})\ge 0\). Similarly, since \(\boldsymbol{u}^{*}\) is a team-optimal solution, it satisfies \(G(\boldsymbol{u}^{*})^\top(\boldsymbol{u}-\boldsymbol{u}^{*})\ge 0\) for all \(\boldsymbol{u}\in \boldsymbol{\Xi}\). By choosing \(\boldsymbol{u}=\boldsymbol{u}^{\Diamond}\), we obtain \(G(\boldsymbol{u}^{*})^\top(\boldsymbol{u}^{\Diamond}-\boldsymbol{u}^{*})\ge 0\). Adding the two inequalities yields
$\bigl(F(\boldsymbol{u}^{\Diamond})-G(\boldsymbol{u}^{*})\bigr)^\top(\boldsymbol{u}^{\Diamond}-\boldsymbol{u}^{*})\le 0$.
Hence,
$\bigl(G(\boldsymbol{u}^{\Diamond})-G(\boldsymbol{u}^{*})\bigr)^\top(\boldsymbol{u}^{\Diamond}-\boldsymbol{u}^{*})
\le
\bigl(G(\boldsymbol{u}^{\Diamond})-F(\boldsymbol{u}^{\Diamond})\bigr)^\top(\boldsymbol{u}^{\Diamond}-\boldsymbol{u}^{*})$.
Since \(G\) is \(\kappa_1\)-strongly monotone, we have
$\kappa_1\|\boldsymbol{u}^{\Diamond}-\boldsymbol{u}^{*}\|^2
\le
\bigl(G(\boldsymbol{u}^{\Diamond})-G(\boldsymbol{u}^{*})\bigr)^\top(\boldsymbol{u}^{\Diamond}-\boldsymbol{u}^{*})$,
and therefore
$\kappa_1\|\boldsymbol{u}^{\Diamond}-\boldsymbol{u}^{*}\|^2
\le
\bigl(G(\boldsymbol{u}^{\Diamond})-F(\boldsymbol{u}^{\Diamond})\bigr)^\top(\boldsymbol{u}^{\Diamond}-\boldsymbol{u}^{*})$.
Applying the Cauchy--Schwarz inequality gives
$\kappa_1\|\boldsymbol{u}^{\Diamond}-\boldsymbol{u}^{*}\|^2
\le
\|F(\boldsymbol{u}^{\Diamond})-G(\boldsymbol{u}^{\Diamond})\|\,\|\boldsymbol{u}^{\Diamond}-\boldsymbol{u}^{*}\|$.
Thus,
\[
\|\boldsymbol{u}^{\Diamond}-\boldsymbol{u}^{*}\|
\le
\frac{1}{\kappa_1}\|F(\boldsymbol{u}^{\Diamond})-G(\boldsymbol{u}^{\Diamond})\|\le
\frac{1}{\kappa_1} \sup_{\boldsymbol{u}\in \boldsymbol{\Xi}}\|F(\boldsymbol{u})-G(\boldsymbol{u})\|.
\]
Next, define for each \(i\in \mathcal{N}\), \(\Delta_i(\boldsymbol{u}) := \nabla_{u_i}\mathcal{C}_i(\boldsymbol{u})-\nabla_{u_i}\mathcal{C}(\boldsymbol{u})\). By the expressions of \(\nabla_{u_i}\mathcal C_i\) and \(\nabla_{u_i}\mathcal C\),
\begin{align*} \|\Delta_i(\boldsymbol{u})\| \le\;& \|\nabla_{u_i}f_i(u_i;\alpha_i)-\nabla_{u_i}f_i(u_i;\alpha)\| +\|\gamma_i-\gamma\| \\ &+\|(\nabla_{u_i}g_i(u_i;\beta_i)-\nabla_{u_i}g_i(u_i;\beta))^\top \sigma(\boldsymbol{u})\| \\ &+\!|\sigma_i|\!\|g_i(u_i;\beta_i)\!-\!g_i(u_i;\beta)\| \!+\!|\sigma_i|\!\sum_{j\neq i}\!\|g_j(u_j;\beta)\|. \end{align*}
Since $\boldsymbol{\Xi}$ is compact, there exist constants $\bar{\sigma}>0$ and $\bar g_j>0$ such that
$\|\sigma(\boldsymbol{u})\|\le \bar{\sigma}$, $\|g_j(u_j;\beta)\|\le \bar g_j$,
$\forall \boldsymbol{u}\in \boldsymbol{\Xi}$,  $j\in \mathcal N$.
Hence, by Assumption~2, $\|\Delta_i(\boldsymbol{u})\|
\le
L_i^{\nabla f}\|\alpha_i-\alpha\|
+
(\bar{\sigma}L_i^{\nabla g}+|\sigma_i|L_i^g)\|\beta_i-\beta\|
+
|\sigma_i|\sum_{j\neq i}\bar g_j
+
\|\gamma_i-\gamma\|$. 
Therefore,
\begin{align*}
\|F(\boldsymbol{u})\!-\!G(\boldsymbol{u})\|
&\!=\!
(\sum_{i=1}^N \!\|\Delta_i(\boldsymbol{u})\|^2)^{1/2} \\
&\!\le\!
(
\sum_{i=1}^N \!
(
L_i^{\nabla f}\!\|\alpha_i\!-\!\alpha\|
\!+\!
(\bar{\sigma}L_i^{\nabla g}\!\!+\!|\sigma_i|L_i^g)\|\beta_i\!-\!\!\beta\|\\
&\quad\quad\!+\!
|\sigma_i|\sum_{j\neq i}\bar g_j
+
\|\gamma_i-\gamma\|
)^2
)^{1/2}
=: \xi .
\end{align*}
Combining the above estimates, for any \(\boldsymbol{u}^{\Diamond}\in \Upsilon_{NE}\), we obtain
$\|\boldsymbol{u}^{\Diamond}-\boldsymbol{u}^{*}\|
\le
\frac{1}{\kappa_1}\sup_{\boldsymbol{u}\in \boldsymbol{\Xi}}\|F(\boldsymbol{u})-G(\boldsymbol{u})\|
\le
\frac{1}{\kappa_1}\,\xi$. 
Since \(G\) is \(\kappa_1\)-strongly monotone, the team-optimal solution is unique, i.e., \(\Upsilon_{TO}=\{\boldsymbol{u}^{*}\}\). Hence,
\[
d_H(\Upsilon_{TO},\Upsilon_{NE})=
\sup_{\boldsymbol{u}^{\Diamond}\in \Upsilon_{NE}}
\|\boldsymbol{u}^{\Diamond}-\boldsymbol{u}^{*}\|
\le
\frac{1}{\kappa_1}\,\xi.
\]
This completes the proof. \hfill $\square$

In Theorem~\ref{t13}, the upper bound $\xi$ on the Hausdorff distance \( d_H(\Upsilon_{TO}, \Upsilon_{NE}) \) depends only on the feasible set and the structures of \(F\) and \(G\). Therefore, it can be estimated directly from the game primitives, thereby providing a tractable worst-case bound on the deviation of the NE set from the team-optimal set without computing all NE and team-optimal solutions. If \(\xi\) is small, then the NE set remains close to the team-optimal set even in the worst case.
 }

The results in this section regarding the consistency and deviation conditions are derived under fixed game parameters. In the following section, we turn to the case where these parameters are adjustable and study how such adjustments can reduce the deviation.




\section{Deviation Reduction Between NE and Team-Optimal Solution}\label{sec:deviation reduction}

When the NE and the team-optimal solution are inconsistent, it is desirable to mitigate their deviation,
as posed in the third question of Section~\ref{sec: formulation}.

To this end, we introduce the notion of a team mediator.  
The mediator’s objective is to reshape the team members’ perceived game parameters to shift the NE closer to team‐optimality. {For example, in  the shortest path traffic problem, the mediator may act as a network operator that modifies vehicles’ perceived local, interaction, or route-dependent travel costs through pricing or incentives, thereby shifting the resulting NE toward the team-optimal flow pattern.}
This naturally leads to a single-leader multi-follower scheme \cite{von2010market}, in which the mediator acts as the leader and the members act as followers. 

Specifically, recall each member's individual parameters in \eqref{para}. For each $i\in\mathcal{N}$, let $\Delta \alpha_i\in\mathbb{R}^{d_\alpha}$, $\Delta \beta_i\in\mathbb{R}^{d_\beta}$, and $\Delta \gamma_i\in \mathbb{R}^{d_\gamma}$ 
denote the parameter adjustments determined by the team mediator. Define  $\theta_i=\operatorname{col}\{\Delta \alpha_i,\Delta \beta_i,\Delta \gamma_i\}\in\mathbb{R}^d$, $d\in\mathbb{N}$,  where $d=d_\alpha+d_\beta+d_\gamma$, and let $\boldsymbol{\theta}=\operatorname{col}\{\theta_i\}_{i=1}^N\in\mathbb{R}^{dN}$ denote the mediator's decision variable.
Suppose that the mediator knows the value of the team-optimal solution $\boldsymbol{u}^*$.
On the upper level, the team mediator manipulates each member's parameters
$\alpha_i$, $\beta_i$, $\gamma_i$ to   $\alpha_i+\Delta \alpha_i$, $\beta_i+\Delta \beta_i$, $\gamma_i+\Delta \gamma_i$, aiming to steer the NE strategy with the team-optimal solution. 
On the lower level, members engage in a manipulated game with  adjusted parameters, in which each member minimizes its manipulated cost function  $\mathcal{C}_i(\theta_i,\boldsymbol{u}): \mathbb{R}^d\times\mathbb{R}^{nN}\rightarrow \mathbb{R}$.

{Incorporating the lower-level solution mapping 
$\boldsymbol{u}^{\Diamond}(\boldsymbol{\theta})$, 
the mediator aims to solve the following optimization problem:
\begin{align}\label{manager_goal}
&\min_{\boldsymbol{\theta}} \;\Psi(\boldsymbol{\theta}):=\frac{1}{2} \|\boldsymbol{u}^{\Diamond}(\boldsymbol{\theta})-\boldsymbol{u}^{*}\|^2
+\frac{\rho}{2}\|\boldsymbol{\theta}\|^2,
\end{align}
where $\rho>0$ is a weight decay parameter for penalizing the quadratic term $\frac{1}{2}\|\boldsymbol{\theta}\|^2$. We denote $\psi(\boldsymbol{\theta}):=\frac{1}{2} \|\boldsymbol{u}^{\Diamond}(\boldsymbol{\theta})-\boldsymbol{u}^{*}\|^2$.
}

\begin{remark}
The objective in \eqref{manager_goal} minimizes the deviation between the NE and the team-optimal solution. This deviation does not necessarily vanish, as the game structure may prevent exact alignment.If the team-optimal solution \(\boldsymbol{u}^*\) is computationally expensive to obtain or otherwise unavailable, one may instead use \(\|F(\boldsymbol{\theta},\boldsymbol{u}^\Diamond(\boldsymbol{\theta}))-G(\boldsymbol{u}^\Diamond(\boldsymbol{\theta}))\|^2\) as an upper-bound surrogate for the equilibrium deviation, since Theorem~\ref{t13} gives \(\|\boldsymbol{u}^\Diamond(\boldsymbol{\theta})-\boldsymbol{u}^*\|\leq \kappa_1^{-1}\|F(\boldsymbol{\theta},\boldsymbol{u}^\Diamond(\boldsymbol{\theta}))-G(\boldsymbol{u}^\Diamond(\boldsymbol{\theta}))\|\). 
\end{remark}


\subsection{Inverse learning}

{

In a leader--follower scheme, the leader is typically assumed to know the followers' cost functions and thus can characterize their best-response mappings. In practice, however, the mediator may have only partial knowledge of the members' cost functions and may observe only their best-response values under given actions, rather than their explicit functional forms. 
This is because the members' weights $\boldsymbol{\sigma}\!=\!(\sigma_1,\!\dots,\! \sigma_N)$ in the aggregative term $\sigma(\boldsymbol{u})\!=\!\sum_{i=1}^N \!\sigma_i u_i$, which characterize each member's contribution to the aggregate effect, may be unknown to the mediator due to environmental uncertainty or intentional concealment \cite{nguyen2009security,birmpas2020optimally}. Therefore, before determining the optimal intervention strategy, the mediator must first infer these unknown aggregation weights so as to construct a lower-level game model that both best explains the observed equilibrium responses and serves as the basis for the subsequent intervention computation.

To this end, we propose a data-driven inverse learning scheme to learn the lower-level game with an unknown weight vector $\boldsymbol{\sigma}$ by exploiting the connection between the observed NE and the associated first-order optimality condition.
 Let \( F(\boldsymbol{\theta}, \boldsymbol{u}) \!:= \!\operatorname{col}\!\left\{ \nabla_{u_i} \mathcal{C}_i(\theta_i, \!\boldsymbol{u}) \right\}_{i=1}^N \) denote the pseudo-gradient mapping of the member game. Recall the feasible set $\Xi_i$ for $i\!\in\!\mathcal{N}$ and the set $\boldsymbol{\Xi}$  for all members. Take $\boldsymbol{\Xi} \!=\! \{\boldsymbol{u}\!\in\!\mathbb R^{nN}\mid D\boldsymbol{u}\!\le\! b, H\boldsymbol{u}\!=\!m\}   $ with $D\!=\!\operatorname{blkdiag}(D_1,\dots,D_N)\in\mathbb{R}^{q\times Nn}$, $H\!=\!\operatorname{blkdiag}(H_1,\dots,H_N)\in\mathbb{R}^{r\times Nn}$, $b\!=\!\operatorname{col}(b_1,\dots,b_N)\in\mathbb{R}^{q}$ and $m\!=\!\operatorname{col}(m_1,\dots,m_N)\in\mathbb{R}^{r}$. Under Assumption~1 (1)-(2), for any given $\boldsymbol{\theta}$, a strategy profile $\boldsymbol{u}^{\Diamond}\in\boldsymbol{\Xi}$ is called an $\epsilon$-NE, with $\epsilon\ge 0$, if
 \begin{align}\label{approne}
      F(\boldsymbol{\theta},\boldsymbol{u}^{\Diamond})^\top (\boldsymbol{u}-\boldsymbol{u}^{\Diamond})\ge -\epsilon,\quad \forall \boldsymbol{u}\in \boldsymbol{\Xi}.
 \end{align}
In particular, when $\epsilon=0$, \eqref{approne} reduces to the standard optimal condition for NE~\cite[Proposition 1.4.2]{facchinei2003finite}.

The following lemma establishes an equivalent inverse relationship derived from \eqref{approne}, which will play an important role later in developing a data-driven learning method.

\begin{lemma}\label{inversevi}
Under Assumption~1 (1), $\boldsymbol{u}^{\Diamond}$ is an $\epsilon$-NE if and only if there exist $\lambda\ge 0_q$ and  $\nu \in \mathbb{R}^r$ such that
\begin{align}\label{eq:dual_poly_2}
D^\top \lambda + H^\top \nu + F(\boldsymbol{\theta},\boldsymbol{u}^{\Diamond}) &= 0_{Nn}, \\
F(\boldsymbol{\theta},\boldsymbol{u}^{\Diamond})^\top \boldsymbol{u}^{\Diamond} + b^\top \lambda + m^\top \nu &\le \epsilon.\notag 
\end{align}
\end{lemma}
\noindent\textbf{Proof.}   
It follows from \eqref{approne} that
$
F(\boldsymbol{\theta},\boldsymbol{u}^{\Diamond})^\top \boldsymbol{u}^{\Diamond}-\epsilon \le F(\boldsymbol{\theta},\boldsymbol{u}^{\Diamond})^\top \boldsymbol{u},\; \forall \boldsymbol{u}\in \boldsymbol{\Xi}.
$
This inequality should be satisfied for all $\boldsymbol{u}\in \boldsymbol{\Xi}$, including
the minimum, namely,
\begin{equation}\label{min}
F(\boldsymbol{\theta},\boldsymbol{u}^{\Diamond})^\top \boldsymbol{u}^{\Diamond}-\epsilon \le \min_{\boldsymbol{u}\in \boldsymbol{\Xi}} F(\boldsymbol{\theta},\boldsymbol{u}^{\Diamond})^\top \boldsymbol{u}.
\end{equation}
Devise a Lagrangian function for the optimization $\min_{\boldsymbol{u}\in \boldsymbol{\Xi}} F(\boldsymbol{\theta},\boldsymbol{u}^{\Diamond})^\top \boldsymbol{u}$ as follows:
\begin{align*}
\mathcal{L}_1(\boldsymbol{u},\lambda,\nu)
&= F(\boldsymbol{\theta},\boldsymbol{u}^{\Diamond})^\top \boldsymbol{u}
+ \lambda^\top(D\boldsymbol{u}-b)
+ \nu^\top(H\boldsymbol{u}-m) \\
&= (F(\boldsymbol{\theta},\boldsymbol{u}^{\Diamond})+D^\top\lambda+H^\top\nu)^\top\boldsymbol{u}-b^\top\lambda-m^\top\nu,
\end{align*}
where $\lambda\in\mathbb{R}^q_+$ and $\nu\in\mathbb{R}^r$ are multipliers.

Under Assumption 1 (1), given any  $\boldsymbol{\theta}$, the duality gap of $\mathcal{L}_1$ vanishes.
Thus, the minimum of $F(\boldsymbol{\theta},\boldsymbol{u}^{\Diamond})^\top \boldsymbol{u}$ subject to $\boldsymbol{u}\in \boldsymbol{\Xi}$ equals the
maximum of the objective for the following constrained optimization:
\[
\max_{\lambda\ge 0_q,\nu} -b^\top\lambda-m^\top\nu
\quad
\text{s.t.}\quad
F(\boldsymbol{\theta},\boldsymbol{u}^{\Diamond})+D^\top\lambda+H^\top\nu=0_{Nn}.
\]
On this basis, the right side in  \eqref{min} can be rewritten as
\begin{equation}\label{dualmax}
F(\boldsymbol{\theta},\boldsymbol{u}^{\Diamond})^\top \boldsymbol{u}^{\Diamond}-\epsilon \le \max_{\lambda,\nu} -b^\top\lambda-m^\top\nu,
\end{equation}
with
\begin{equation}\label{dualconst}
F(\boldsymbol{\theta},\boldsymbol{u}^{\Diamond})+D^\top\lambda+H^\top\nu=0_{nN},\;\lambda\ge 0_q.
\end{equation}
If the inequality in \eqref{dualmax} needs to be satisfied with the right-hand side taking on
the maximum, then we should merely ensure that there exists at least one feasible
$\lambda$ and $\nu$ satisfy \eqref{dualconst}. The inequality in \eqref{dualmax} finally becomes
\begin{align*}
   & F(\boldsymbol{\theta},\boldsymbol{u}^{\Diamond})^\top \boldsymbol{u}^{\Diamond}-\epsilon \le  -b^\top\lambda-m^\top\nu, \\
   & F(\boldsymbol{\theta},\boldsymbol{u}^{\Diamond})+D^\top\lambda+H^\top\nu=0_{Nn},\;\lambda\ge 0_q \notag.
\end{align*}
The reverse proof  can be conducted similarly using
weak duality properties considering the aforementioned equivalent conversions.
 \hfill $\square$

Based on Lemma~\ref{inversevi}, we are now in a position to design a data-driven inverse learning problem. Recalling the unknown parameter $\boldsymbol{\sigma}$ in the pseudo-gradient mapping $F$ in \eqref{gradientne}, we make its parametric dependence explicit and write the lower-level mapping as $F(\boldsymbol{\theta},\boldsymbol{u};\boldsymbol{\sigma})$. Moreover, we assume that $F(\boldsymbol{\theta},\boldsymbol{u};\boldsymbol{\sigma})$ depends continuously on $\boldsymbol{\sigma}$. The goal of the inverse-learning stage is to construct a behaviorally consistent lower-level game model from observed equilibrium samples. In particular, Lemma~\ref{inversevi} converts the original equilibrium condition, which is a semi-infinite constraint over all $\boldsymbol{u}\in\boldsymbol{\Xi}$, into a finite-dimensional certificate involving the dual variables $\lambda$ and $\nu$. This reformulation makes it possible to infer a lower-level game model directly from observed equilibrium responses.

Accordingly, the mediator seeks a parameter vector $\boldsymbol{\sigma}$ such that the resulting game model best explains the observed equilibrium responses with the smallest overall inconsistency. To this end, the mediator may deliberately take a sequence of different actions, for instance in a randomized manner, to generate a collection of observations $(\boldsymbol{\theta}[j],\boldsymbol{u}^{\Diamond}[j])$ for $j\in\mathcal{M}=\{1,\dots,M\}$. Based on these actions and the associated equilibrium responses, the mediator solves
\begin{align}\label{data-driven}
&\min_{\boldsymbol{\epsilon},\boldsymbol{\lambda},\boldsymbol{\nu},\boldsymbol{\sigma}} \|\boldsymbol{\epsilon}\|\\
&\operatorname{s.t.}\;
D^\top \lambda[j] + H^\top \nu[j] + F(\boldsymbol{\theta}[j], \boldsymbol{u}^{\Diamond}[j];\boldsymbol{\sigma}) = 0_{Nn}, \quad j \in \mathcal{M}, \notag\\
&\quad\;\;\;\;
F(\boldsymbol{\theta}[j],\! \boldsymbol{u}^{\Diamond}[j];\!\boldsymbol{\sigma})^\top \!\boldsymbol{u}^{\Diamond}[j]
\!+\! b^\top \lambda[j] \!+ \!m^\top \nu[j] \!\le\! \epsilon[j], \; j\!\in\!\mathcal{M}, \notag\\
&\quad\;\;\;\;
\epsilon[j]\ge 0,\;\lambda[j] \ge 0_q,\quad j\in\mathcal{M}. \notag
\end{align}
Here, $\boldsymbol{\epsilon}=\operatorname{col}\{\epsilon[j]\}_{j=1}^M$, $\boldsymbol{\lambda}=\operatorname{col}\{\lambda[j]\}_{j=1}^M$, and $\boldsymbol{\nu}=\operatorname{col}\{\nu[j]\}_{j=1}^M$, and $\|\cdot\|$ denotes a chosen norm. In problem~\eqref{data-driven}, $\boldsymbol{\epsilon}$ is optimized jointly with the unknown $\boldsymbol{\sigma}$ and measures the degree to which the observed decisions fail to satisfy the exact equilibrium conditions under the candidate model. Hence, minimizing $\|\boldsymbol{\epsilon}\|$ seeks a learned parameter $\hat{\boldsymbol{\sigma}}$ such that the induced lower-level game model best explains the observations with the smallest overall inconsistency. The learned parameter $\hat{\boldsymbol{\sigma}}$ is then used to construct the estimated response map for the subsequent upper-level intervention.

 Note that $\boldsymbol{u}^{\Diamond}$ in \eqref{data-driven} are treated as known data rather than decision variables. Therefore, the complexity of problem \eqref{data-driven} depends on the dependence of $F$ on $\boldsymbol{\sigma}$, rather than on its dependence on $\boldsymbol{u}$. Since $F$ is linear in $\boldsymbol{\sigma}$ in our setting, problem \eqref{data-driven} is convex and hence computationally tractable.

 
\begin{remark}

The learned game model obtained from~\eqref{data-driven} is used not only to explain the observed NE samples, but also to provide a behaviorally consistent estimated response map for the subsequent intervention design. On the one hand, by solving~\eqref{data-driven}, the mediator learns a lower-level game model that fits the observed equilibrium samples by minimizing the overall inconsistency measured by $\boldsymbol{\epsilon}$, which in turn indicates that the equilibrium induced by the learned model is close to the observed one. On the other hand, the learned model enjoys generalization guarantees~\cite{bertsimas2015data,chen2025inverse}, in the sense that it continues to approximately satisfy the equilibrium conditions on unseen data, typically measured by equilibrium inconsistency or violation probability; see~\cite{bertsimas2015data,chen2025inverse} for further details. In our setting, this predictive interpretation is particularly relevant because the subsequent upper-level intervention may involve values of $\boldsymbol{\theta}$ that are not used in the inverse-learning stage. Hence, the learned equilibrium map serves as a meaningful surrogate response map for the subsequent upper-level intervention.
\end{remark}

}
\vspace{-0.2cm}
\subsection{Adam-type hypergradient update}

{
After solving inverse learning problem~\eqref{data-driven}, we take 
 $\hat{\boldsymbol{\sigma}}=(\hat{\sigma}_1,\dots,\hat{\sigma}_N)$ as the optimal solution of \eqref{data-driven} with $\hat{\sigma}(\boldsymbol{u})=\sum_{i=1}^N \hat{\sigma}_i u_i$ as the expression of the estimated aggregator. Then, 
the mediator aims to solve the following problem
\begin{align}\label{manager_goal_est}
&\min_{\boldsymbol{\theta}} \;\Psi(\boldsymbol{\theta};\hat{\boldsymbol{\sigma}}):=\frac{1}{2} \|\boldsymbol{u}^{\Diamond}(\boldsymbol{\theta};\hat{\boldsymbol{\sigma}})-\boldsymbol{u}^{*}\|^2
+\frac{\rho}{2}\|\boldsymbol{\theta}\|^2.
\end{align}}
We first observe that in certain special cases, the mediator can derive a closed-form expression for the optimal parameter adjustments. 
 Specifically,  by ensuring that each member’s cost function 
satisfies condition \eqref{potential}, the induced NE coincides with the team-optimal solution, as established in Proposition~\ref{p1}. The following example illustrates this alignment.

{\noindent\textbf{Example of closed-form solutions}\;
Consider a class of linear-quadratic team and game problems, where the team objective and individual objectives are given by:
\begin{align}\label{para_quad}
\vspace{-0.2cm}
&\mathcal C(\boldsymbol u)
=
\sum_{i=1}^N
(
u_i^\top Q(\alpha)u_i
+
u_i^\top Q(\beta)\hat{\sigma}(\boldsymbol u)
+
\gamma^\top u_i),
\notag\\
   & \mathcal C_i(\theta_i,\boldsymbol u;\hat{\boldsymbol{\sigma}})
=
u_i^\top Q(\alpha_i+\Delta\alpha_i)u_i
+
u_i^\top Q(\beta_i+\Delta\beta_i)\hat{\sigma}(\boldsymbol u)\notag\\
&\quad\quad\quad\quad\quad\quad+
(\gamma_i+\Delta\gamma_i)^\top u_i,
    \end{align}
where $\hat{\sigma}(\boldsymbol{u})\!=\!\frac{1}{N}\sum_{i=1}^N\! u_i$,  and $\alpha,\beta\!\in\!\mathbb R^{d_Q}$ are parameterized over the same linearly independent basis $\{Q_l\}_{l=1}^{d_Q}\!\!\subset\! \mathbb S_{++}^n$, i.e.,
$Q(x)\!\!=\!\sum_{l=1}^{d_Q}\!\! Q_l [x]_l$,  $x\in\mathbb R^{d_Q}$.
By calculating  $ \nabla_{u_i}\!\mathcal C_i\!=\!
2Q(\alpha_i\!+\!\Delta\alpha_i)u_i
\!+\!
Q(\beta_i\!+\!\Delta\beta_i)\hat{\sigma}(\boldsymbol u)
\!+\!
\frac{1}{N}Q(\beta_i\!+\!\Delta\beta_i)u_i
\!+\!
\gamma_i\!+\!\Delta\gamma_i$ and $\nabla_{u_i}\!\mathcal C\!
\!=\!
2Q(\alpha)u_i\!+\!2Q(\beta)\hat{\sigma}(\boldsymbol u)\!+\!\gamma$, one explicit form of the optimal parameter adjustment is $\boldsymbol{\theta}\!=\!\operatorname{col}\{\operatorname{col}\{\alpha\!-\!\tfrac{1}{N}\beta\!-\!\alpha_i,\;2\beta\!-\!\beta_i,\;\gamma\!-\!\gamma_i\}\}_{i=1}^N\!$
for $i \!\in\!\mathcal{N}$ such that  $\nabla_{u_i}\mathcal{C}_i
\!=\!
\nabla_{u_i}\mathcal C$. \hfill$\square$ }%

{However,  in general settings,
a closed-form solution of \eqref{manager_goal_est} is difficult to obtain.  This is due to the nonlinear coupling among the cost parameters and the fact that $\boldsymbol{u}^{\Diamond}(\boldsymbol{\theta};\hat{\boldsymbol{\sigma}})$ typically lacks an explicit form and may be non-differentiable, rendering the composite function $\Psi(\boldsymbol{\theta};\hat{\boldsymbol{\sigma}})$ non-convex and non-smooth. To this end, we adopt an Adam optimizer to obtain the mediator's strategy, in which the key idea is to compute the gradient of $\Psi(\boldsymbol{\theta};\hat{\boldsymbol{\sigma}})$, known as the hypergradient. }
If $\boldsymbol{u}^{\Diamond}(\cdot)$ is differentiable at $\boldsymbol{\theta}$, then by the chain rule, the gradient $\nabla\Psi(\boldsymbol{\theta};\hat{\boldsymbol{\sigma}})$ is given by
\begin{align*}
    \nabla\Psi(\boldsymbol{\theta};\hat{\boldsymbol{\sigma}})\!= \!\nabla\psi(\boldsymbol{\theta};\hat{\boldsymbol{\sigma}})\!+\!\rho\boldsymbol{\theta}\!=\!\text{J}{\boldsymbol{u}^{\Diamond}}\!(\boldsymbol{\theta};\hat{\boldsymbol{\sigma}})^\top \!(\boldsymbol{u}^{\Diamond}(\boldsymbol{\theta};\hat{\boldsymbol{\sigma}})\!\!-\!\!\boldsymbol{u}^{*})\!+\!\rho\boldsymbol{\theta},
\end{align*}
where $\text{J}{\boldsymbol{u}^{\Diamond}}(\boldsymbol{\theta};\hat{\boldsymbol{\sigma}})$ denotes Jacobian of $\boldsymbol{u}^{\Diamond}$ with respect to $\boldsymbol{\theta}$.

Computing $\nabla\Psi(\boldsymbol{\theta};\hat{\boldsymbol{\sigma}})$ 
requires evaluating $\boldsymbol{u}^{\Diamond}(\boldsymbol{\theta};\hat{\boldsymbol{\sigma}})$ and  its Jacobian $\text{J}{\boldsymbol{u}^{\Diamond}}(\boldsymbol{\theta};\hat{\boldsymbol{\sigma}})$. To obtain the NE $\boldsymbol{u}^{\Diamond}(\boldsymbol{\theta};\hat{\boldsymbol{\sigma}})$, we employ a projected gradient method:
\begin{equation}\label{teamNE}	
	{u}_i[l+1]=\Pi_{{\Xi}_i}({u}_i[l]-\tau \nabla_{u_i}\mathcal{C}_i({{\theta}}_i,\boldsymbol{u}[l])),\; \forall l\in \mathbb{N},	
		\end{equation}
    with an  initial $ {u}_i[0]\in \mathbb{R}^n$,  where $\tau>0$   is the stepsize. In compact form, \eqref{teamNE}  becomes
     $\boldsymbol{u}[l+1]=\Pi_{\boldsymbol{\Xi}}(\boldsymbol{u}[l]-\tau F(\boldsymbol{\theta},\boldsymbol{u}[l]))$.		
         Under Assumption \ref{assum2},  $\boldsymbol{u}^{\Diamond}(\boldsymbol{\theta};\hat{\boldsymbol{\sigma}})$ is equivalently characterized as a
fixed point of the projected gradient map: 
\begin{equation}\label{fixed}
    \boldsymbol{u}^{\Diamond}\!(\boldsymbol{\theta};\!\hat{\boldsymbol{\sigma}})\!=\!\Pi_{\boldsymbol{\Xi}}\!(\boldsymbol{u}^{\Diamond}\!(\boldsymbol{\theta};\hat{\boldsymbol{\sigma}})\!-\!\tau F\!(\boldsymbol{\theta},\!\boldsymbol{u}^{\Diamond}(\boldsymbol{\theta};\!\hat{\boldsymbol{\sigma}})))\!:=\!\zeta\!(\boldsymbol{\theta},\!\boldsymbol{u}^{\Diamond}(\boldsymbol{\theta};\hat{\boldsymbol{\sigma}})).
\end{equation}
    To evaluate the Jacobian $\text{J}{\boldsymbol{u}^{\Diamond}}(\boldsymbol{\theta};\hat{\boldsymbol{\sigma}})  $, we are intuitively motivated by differentiating \eqref{fixed} and adopt a fixed-point iterative scheme of the form: 
\begin{equation}\label{jaco}	
z[p+1]=\text{J}_{\boldsymbol{u}^{\Diamond}}\zeta(\boldsymbol{\theta},\boldsymbol{u}^{\Diamond})z[p]+\text{J}_{\boldsymbol{\theta}}\zeta(\boldsymbol{\theta},\boldsymbol{u}^{\Diamond}),\; \forall p\in \mathbb{N},
		\end{equation}
with an  initial $ z[0]\!\in\! \mathbb{R}^{nN\times dN}$, where  all terms are evaluated at $\boldsymbol{u}^{\Diamond}\!\!\!:=\!\boldsymbol{u}^{\Diamond}(\boldsymbol{\theta};\!\hat{\boldsymbol{\sigma}})$. The Jacobians in \eqref{jaco} are given by
%
   \begin{align*}
\text{J}_{\boldsymbol{u}^{\Diamond}}\zeta(\boldsymbol{\theta},\boldsymbol{u}^{\Diamond})\!& =\! \text{J}\phi(\boldsymbol{u}^{\Diamond}-\tau F(\boldsymbol{\theta},\boldsymbol{u}^{\Diamond};\!\hat{\boldsymbol{\sigma}})) \left( I \!-\! \tau \text{J}_{\boldsymbol{u}^{\Diamond}}F(\boldsymbol{\theta}, \boldsymbol{u}^{\Diamond};\!\hat{\boldsymbol{\sigma}}) \right),\\
\text{J}_{\boldsymbol{\theta}}\zeta(\boldsymbol{\theta},\boldsymbol{u}^{\Diamond}) &= - \tau \text{J}\phi(\boldsymbol{u}^{\Diamond}-\tau F(\boldsymbol{\theta},\boldsymbol{u}^{\Diamond};\!\hat{\boldsymbol{\sigma}}))\text{J}_{\boldsymbol{\theta}}F(\boldsymbol{\theta}, \boldsymbol{u}^{\Diamond};\!\hat{\boldsymbol{\sigma}}),
\end{align*}
 where $\text{J}\phi(\cdot)$ denotes the Jacobian of the projection operator $\phi(\cdot):=\Pi_{\boldsymbol{\Xi}}(\cdot)$. Since $\boldsymbol{\Xi}$ is a polyhedron,  $\phi(\cdot)$ can be characterized as the solution to a standard quadratic program.  The evaluation of $\text{J}\phi(\cdot)$ reduces to solving a linear system of equations by differentiating the Karush-Kuhn-Tucker (KKT)  conditions of this quadratic program at the solution point, and the resulting linear system is always solvable \cite{amos2017optnet}.
 Analytical characterization and computational methods for computing $\text{J}\phi(\cdot)$ can be found in  \cite{amos2017optnet}. 

 In the case where $\boldsymbol{u}^{\Diamond}(\cdot;\hat{\boldsymbol{\sigma}})$ is non-differentiable at $\boldsymbol{\theta}$, the analysis is similar. 
 We briefly note that both $\boldsymbol{u}^{\Diamond}(\cdot;\hat{\boldsymbol{\sigma}})$ and $\psi$ are definable, as will be shown in Appendix~\ref{proofthm3}. 
This property ensures that $\psi$ and $\Psi$ are path-differentiable and admit  conservative gradients $\mathcal{J}_{\psi}(\boldsymbol{\theta};\hat{\boldsymbol{\sigma}})$ and $\mathcal{J}_{\Psi}(\boldsymbol{\theta};\hat{\boldsymbol{\sigma}})$, 
to which a generalized chain rule still applies \cite{bolte2021nonsmooth}.
Moreover, by \cite[Theorem 1]{bolte2021conservative}, $\mathcal{J}_{\Psi}(\boldsymbol{\theta};\hat{\boldsymbol{\sigma}}) = \{\nabla \Psi(\boldsymbol{\theta};\hat{\boldsymbol{\sigma}})\}$ almost everywhere, and  in particular, $\mathcal{J}_{\Psi}(\boldsymbol{\theta};\hat{\boldsymbol{\sigma}})=\{\nabla \Psi(\boldsymbol{\theta};\hat{\boldsymbol{\sigma}})\}$ whenever $\Psi$ is differentiable at $\boldsymbol{\theta}$.
We further elaborate on this extension in Appendix~\ref{proofthm3}.

{ On this basis, to improve computational efficiency and numerical stability, we let the mediator update its decision variable $\boldsymbol{\theta}$ using an Adam-type update as follows:
\begin{equation}\label{adam}
\left\{
\begin{aligned}
&\omega[k] \in \mathcal{J}_{\psi}(\boldsymbol{\theta}[k]), 
\\
&\mu[k+1] = (1 - \vartheta_1[k])\,\mu[k] + \vartheta_1[k]\,\omega[k],  \\
&s[k+1] = (1 - \vartheta_2)s[k] + \vartheta_2 \omega^{2}[k], \\
&\boldsymbol{\theta}[k+1] = \boldsymbol{\theta}[k] \!-\! \eta(\sqrt{s[k\!+\!1]} + \varpi)^{-1}\!\odot\!\bigl(\mu[k+1] \!+\! \rho\,\boldsymbol{\theta}[k]\bigr), 
\end{aligned}
\right.
\end{equation}
where  initial point $\boldsymbol{\theta}[0]\in\mathbb{R}^{dN}$, $\mu[0]\in\mathbb{R}^{dN}$, $s[0]\in\mathbb{R}^{dN}_{+}$,  $\rho>0$ is the weight‐decay parameter, $\varpi>0$ is the  safeguard parameter. The operator $\odot$ denotes element-wise multiplication, $\omega^2[k]$ denotes
the element-wise square of $\omega[k]$, and the inverse $(\sqrt{s[k+1]}+\varpi)^{-1}$ is taken componentwise. The sequences \({\mu[k]}\) and \({s[k]}\) are the first- and second-moment estimators of the hypergradient, respectively: \(\mu[k+1]\) smooths the descent direction by tracking an exponential moving average of the hypergradient, while \((\sqrt{s[k+1]}+\varpi)^{-1}\) adaptively rescales each coordinate using the accumulated squared hypergradient. The sequences   $\{\vartheta_1[k]\}_{k\in\mathbb{N}}$, $\vartheta_2\in(0,1)$ and $\eta$ are the stepsizes  for the momentum term $\{\mu[k]\}$, the estimator $\{s[k]\}$ and   the decision variable $\{\boldsymbol{\theta}[k]\}$.

\begin{algorithm}[t]
	\caption{}
	\label{centralized-alg}
	\begin{algorithmic}[1]
  \renewcommand{\algorithmicrequire}{ 
  \textbf{Initialize:}}
		\REQUIRE Set $\boldsymbol{\theta}[0]\in\mathbb{R}^{dN}$, $\boldsymbol{u}[0]\in\mathbb{R}^{nN}$, $z[0]\in \mathbb{R}^{nN\times dN}$, $\mu[0]\in\mathbb{R}^{dN}$, $s[0]\in\mathbb{R}^{dN}_{+}$, $ \varpi>0$, $\vartheta_2\in(0,1)$, $\rho>0$
       \FOR{$j = 1,2,\dots M$}
       \STATE Collect observations $(\boldsymbol{\theta}[j],\boldsymbol{u}^{\Diamond}(\boldsymbol{\theta}[j]))$
        \ENDFOR
        \STATE Solve inverse problem \eqref{data-driven} and obtain $\hat{\boldsymbol{\sigma}}$
		\FOR{$k = 1,2,\dots$}
        \STATE Set ${\boldsymbol{\theta}}=\boldsymbol{\theta}[k]$
         \STATE  \textbf{Inner loop}:
\quad Compute the NE  ${	\boldsymbol{u}^{\Diamond}({\boldsymbol{\theta}};\hat{\boldsymbol{\sigma}})}$

         \STATE 
         \textbf{Outer loop:}

      \quad Evaluate the Jacobian $\text{J}{\boldsymbol{u}^{\Diamond}}(\boldsymbol{\theta};\hat{\boldsymbol{\sigma}})$
    
   \quad Update $\boldsymbol{\theta}[k+1]$ via ${\boldsymbol{u}^{\Diamond}({\boldsymbol{\theta}};\hat{\boldsymbol{\sigma}})}$ and $\text{J}{\boldsymbol{u}^{\Diamond}}(\boldsymbol{\theta};\hat{\boldsymbol{\sigma}})$
	\ENDFOR
	\end{algorithmic}
 
\end{algorithm}

For clarity, we summarize the proposed scheme in Algorithm~1. The mediator first takes a sequence of random actions to collect the members' equilibria, and then solves the inverse problem \eqref{data-driven} to 
obtain the parameter estimate $\hat{\boldsymbol{\sigma}}$. Then the mediator  updates its action using the Adam optimizer~\eqref{adam}, where the hypergradient is constructed based on the corresponding NE strategies and the associated Jacobian.}

\vspace{-0.2cm}

\subsection{Convergence analysis}

In this section, we  study the convergence of \eqref{adam}. 

\begin{assumption}\label{assum3}
\leavevmode
\begin{enumerate}[(1)]
   \item For any fixed \( \boldsymbol{\theta} \), the pseudo-gradient mapping \( F(\boldsymbol{\theta}, \boldsymbol{u};\hat{\boldsymbol{\sigma}}) \) is \( \kappa_2 \)-strongly monotone and \( \nu_2 \)-Lipschitz continuous in \( \boldsymbol{u} \).
    
    \item The mapping \( F(\boldsymbol{\theta}, \boldsymbol{u};\hat{\boldsymbol{\sigma}}) \) is definable\footnote{Definable functions form a broad class that includes most functions used in optimization and machine learning, such as semialgebraic functions, as well as functions involving exponentials and logarithms. Definable functions are closed under standard operations (e.g., addition, multiplication, composition) and possess desirable properties such as path differentiability \cite{bolte2021nonsmooth}.}, continuously differentiable, and there exist constants \( \nu_{\theta} > 0 \), \( \nu_{u} > 0 \) such that the partial Jacobians of \( F \) satisfy:
    \begin{align*}
    &\| \mathrm{J}_{\boldsymbol{\theta}} F(\boldsymbol{\theta}, \boldsymbol{u};\hat{\boldsymbol{\sigma}}) - \mathrm{J}_{\boldsymbol{\theta}} F(\boldsymbol{\theta}, \boldsymbol{u}';\hat{\boldsymbol{\sigma}}) \| \leq \nu_{\theta} \| \boldsymbol{u} - \boldsymbol{u}' \|,\\
    &
    \| \mathrm{J}_{\boldsymbol{u}} F(\boldsymbol{\theta}, \boldsymbol{u};\hat{\boldsymbol{\sigma}}) - \mathrm{J}_{\boldsymbol{u}} F(\boldsymbol{\theta}, \boldsymbol{u}';\hat{\boldsymbol{\sigma}}) \| \leq \nu_{u} \| \boldsymbol{u} - \boldsymbol{u}' \|,
    \end{align*}
    for all \( \boldsymbol{u}, \boldsymbol{u}' \in \boldsymbol{\Xi} \).
    \item {$\psi$ is bounded from below. }

    \item   {There exists a constant $G>0$ such that $\|\boldsymbol{\omega}\|\le G$ for all $\boldsymbol{\theta}\in\mathbb{R}^{dN}$ and all $\boldsymbol{\omega}\in \mathcal{J}_{\psi}(\boldsymbol{\theta};\hat{\boldsymbol{\sigma}})$.}
\end{enumerate}
\end{assumption}

Assumption~\ref{assum3}~(1), together with Assumption~1~(1)-(2), guarantees the existence and uniqueness of an NE for any fixed $\boldsymbol{\theta}$ under the learned parameter $\hat{\boldsymbol{\sigma}}$~\cite[Theorem 2.3.3(b)]{facchinei2003finite}, yielding a single-valued solution mapping $\boldsymbol{u}^{\Diamond}(\boldsymbol{\theta};\hat{\boldsymbol{\sigma}}) : \mathbb{R}^{dN} \to \mathbb{R}^{nN}$. We also restrict our attention to the case where the team-optimal solution is unique in the section. 
{Also, the structural conditions imposed on \( F(\boldsymbol{\theta}, \boldsymbol{u}) \) in Assumption~\ref{assum3}\,(2) are satisfied by a broad class of models~\cite{grontas2024big,bolte2021nonsmooth}, including, for instance, the linear-quadratic cost function introduced in the preceding example. } {Assumption~\ref{assum3}\,(4) implies that the hypergradient remains uniformly bounded over all iterates.}

As for the inner loop,
the convergence of  $\eqref{teamNE}$  is a classical result under
a standard stepsize condition \cite{bauschke2017convex}. We include this result in the following lemma for completeness.
\begin{lemma}\label{con_ne}
Given Assumptions~ 1 (1)-(2) and \ref{assum3}, let the stepsize \( \tau \) satisfy \( 0 < \tau < \frac{2\kappa_2}{\nu_2^2} \). Then, the sequence \( \{ \boldsymbol{u}[l] \}_{l \in \mathbb{N}} \) generated by~\eqref{teamNE} converges linearly to \( \boldsymbol{u}^\Diamond(\boldsymbol{\theta};\hat{\boldsymbol{\sigma}}) \), with a rate given by
$
\sqrt{1 - \tau (2\kappa_2 - \tau \nu_2^2)}.
$
\end{lemma}



For the outer loop, we first study the convergence of the Jacobian iteration in (\ref{jaco}). Note that the behavior of  \eqref{jaco} is related to the differentiability of  $\boldsymbol{u}^{\Diamond}(\boldsymbol{\theta};\hat{\boldsymbol{\sigma}})$. Recall that 
$\boldsymbol{u}^{\Diamond}(\boldsymbol{\theta};\hat{\boldsymbol{\sigma}})=\!\Pi_{\boldsymbol{\Xi}}(\boldsymbol{u}^{\Diamond}(\boldsymbol{\theta};\hat{\boldsymbol{\sigma}})\!-\!\tau F(\boldsymbol{\theta},\boldsymbol{u}^{\Diamond}(\boldsymbol{\theta};\hat{\boldsymbol{\sigma}})))$
in \eqref{fixed}.
Given that $F$ is continuously differentiable,  we only need to consider the differentiability of projection operator $\phi(\cdot)=\Pi_{\boldsymbol{\Xi}}(\cdot)$.  
Note that $\phi(\cdot)$ can be characterized as the solution to a standard quadratic program,
and is therefore  piecewise affine \cite{tondel2003further,amos2017optnet}. That is, $\phi(\cdot)$ is affine on each polyhedral region of its domain and hence is continuously differentiable on the interior of each region. Non-differentiability may only occur on the boundaries between regions, which form a lower-dimensional subspace and are  Lebesgue-negligible.

To proceed, take $\varsigma(\boldsymbol{\theta},\boldsymbol{u};\hat{\boldsymbol{\sigma}}):=\boldsymbol{u}-\tau F(\boldsymbol{\theta},\boldsymbol{u};\hat{\boldsymbol{\sigma}})$. 
The following lemma derives a condition
for the differentiability of $\boldsymbol{u}^{\Diamond}(\boldsymbol{\theta};\hat{\boldsymbol{\sigma}})$ and establishes the convergence of (\ref{jaco}). The proof is given in Appendix \ref{jaco_con}.

\begin{lemma}\label{l11}
Under Assumptions~\ref{assum2} (1)-(2) and~\ref{assum3}, fix any \( \boldsymbol{\theta}  \), and suppose the projection operator \( \phi(\cdot) \) is continuously differentiable at \( \varsigma(\boldsymbol{\theta}, \!\boldsymbol{u}^\Diamond(\boldsymbol{\theta};\!\hat{\boldsymbol{\sigma}})) \). Then, the mapping \( \boldsymbol{u}^\Diamond(\boldsymbol{\theta};\hat{\boldsymbol{\sigma}}) \) is continuously differentiable at \( \boldsymbol{\theta} \). Moreover, the sequence \( \{ z[p] \}_{p \in \mathbb{N}} \) generated by~\eqref{jaco} converges to the Jacobian \( \mathrm{J} \boldsymbol{u}^\Diamond\!(\boldsymbol{\theta};\!\hat{\boldsymbol{\sigma}}) \).
\end{lemma}

\begin{remark}
In the algorithm implementation, we observe that $\phi(\cdot)$ is continuously differentiable at $\varsigma(\boldsymbol{\theta}[k],\boldsymbol{u}^{\Diamond}(\boldsymbol{\theta}[k]))$ 
for all but finitely many iterates ${\boldsymbol{\theta}[k]}$.
In the event that $\phi(\cdot)$ is non-differentiable at a given iterate, the Jacobian in the recursion  \eqref{jaco} is replaced by an arbitrary selection from its conservative Jacobian,
and the fixed-point scheme converges to an element of the conservative Jacobian of $\boldsymbol{u}^{\Diamond}(\boldsymbol{\theta};\hat{\boldsymbol{\sigma}})$.
A formal discussion of this case and the convergence analysis can be found in \cite[Proposition 3]{grontas2024big}.
\end{remark}

{
Next, we consider the convergence of iteration \eqref{adam} for $\boldsymbol{\theta}$. 
In general, first-order methods in non-smooth and non-convex settings are not guaranteed to converge to local minima. Instead, we show that \eqref{adam} converges to critical points of the problem \eqref{manager_goal}, i.e.,
any $\boldsymbol{\theta}$ that satisfies the inclusion
\begin{equation}\label{inclusion1}
\boldsymbol{0}_{dN}\in \mathcal{J}_{\Psi}(\boldsymbol{\theta};\hat{\boldsymbol{\sigma}}).
\end{equation}
Condition \eqref{inclusion1} is necessary for $\boldsymbol{\theta}$ to be a local minimum. 

Let $\operatorname{crit}_{\Psi}= \left\{ {\theta} \in \mathbb{R}^{dN} \;:\; \boldsymbol{0}_{dN} \in \mathcal{J}_{\Psi}(\boldsymbol{\theta};\hat{\boldsymbol{\sigma}})\right\}$ denote
 the set of $\mathcal{J}_{\Psi}({\theta};\hat{\boldsymbol{\sigma}})$-critical points. Now, we are ready to prove the convergence of \eqref{adam}.

\begin{mythm}\label{t5}
Given Assumptions~1 (1)-(2) and~\ref{assum3}, let  the momentum sequence $\{\vartheta_1[k]\}$ satisfy
$\sum_{k=0}^{\infty}\vartheta_1[k]=\infty$, $
\vartheta_1[k]\log k \to 0 $ and the stepsize $\eta$ satisfies $0<\eta<\frac{\varpi}{\rho}$. Then,
every limit point of $\{\boldsymbol{\theta}[k]\}$ is a $\mathcal{J}_{\Psi}({\theta};\hat{\boldsymbol{\sigma}})$-critical point of $\Psi$ and  the sequence $\{\Psi(\boldsymbol{\theta}[k];\hat{\boldsymbol{\sigma}})\}$ converges.
\end{mythm}

We next consider the case where the momentum parameter \(\vartheta_1[k]\) is chosen as a non-diminishing sequence. In this setting, although asymptotic convergence to the critical set may not be guaranteed, we can still establish a neighborhood result.

\begin{corollary}\label{cor:const_momentum} Assume that \(0<\eta<\varpi/\rho\), and let \(\vartheta_1[k]\equiv\bar\vartheta\) for all \(k\ge0\). Then, given Assumptions~1 (1)-(2) and~\ref{assum3}, for any \(\varepsilon>0\), there exists \(\hat\vartheta(\varepsilon)\in(0,1)\) such that, whenever \(0<\bar\vartheta\le \hat\vartheta(\varepsilon)\), the sequence \(\{\theta[k]\}_{k\in\mathbb N}\) generated by \eqref{adam} satisfies \(\limsup_{k\to\infty}\operatorname{dist}(\theta[k],\operatorname{crit}_{\Psi})\le \varepsilon\).
\end{corollary} }

{Finally, we further show that the iteration \eqref{adam} converges to global minima in certain specific settings.}

\noindent\textbf{Example of globally optimal adjustment}\;
Recall the linear-quadratic setting in \eqref{para_quad} and suppose that 
only the  parameters $\operatorname{col}\{ \gamma_i\}_{i=1}^N$ are adjustable.
In this case,  the mediator's decision variable reduces to $\boldsymbol{\theta}=\operatorname{col}\{\Delta \gamma_i\}_{i=1}^N$. {Consider that each strategy set is defined only by affine constraints, i.e., $\Xi_i=\{u_i\in\mathbb{R}^{n}|H_i u_i=m_i\}$. }
Then, solving the KKT conditions of the manipulated game shows that the NE  depends affinely on $\boldsymbol{\theta}$ and takes the form $\boldsymbol{u}^{\Diamond}(\boldsymbol{\theta};\hat{\boldsymbol{\sigma}})=P\boldsymbol{\theta}+p$, where 
$P\in\mathbb{R}^{nN\times dN}$ and $p\in\mathbb{R}^{nN}$ are constants. 
Consequently, the mediator’s objective function is convex in 
$\boldsymbol{\theta}$ and takes the form  $\Psi(\boldsymbol{\theta};\hat{\boldsymbol{\sigma}})=\frac{1}{2}\|P\boldsymbol{\theta}+p-\boldsymbol{u}^*\|^2+\frac{\rho}{2}\|\boldsymbol{\theta}\|^2$. 
This ensures that the mediator achieves a globally optimal adjustment.

{\section{Discussion}\label{sec:diss}
In this section, we provide further discussion on the team-optimal solution, the computational complexity of Algorithm~1 and the comparison with other bi-level methods.  

The team-optimal solution characterizes the optimal joint action that minimizes the team cost $\mathcal{C}$. A related but distinct concept is the social optimum~\cite{koutsoupias1999worst}, which minimizes $\sum_{i=1}^N \mathcal{C}_i$. In general, these two notions do not coincide, since $\mathcal{C}$ may encode a system-level objective and need not equal $\sum_{i=1}^N \mathcal{C}_i$. In our setting, although the team and game costs have similar structural forms, they differ because the team-level coefficients generally satisfy $\alpha\neq \alpha_i$, $\beta\neq \beta_i$, and $\gamma\neq \gamma_i$. The social-welfare cost is recovered only as a special case when the team objective coincides with the sum of the members' local objectives. The inefficiency of NE strategies relative to the social optimum has been extensively studied through the concept of the price of anarchy~\cite{koutsoupias1999worst}, whereas this paper focuses on their inefficiency relative to the team optimum.



The computational complexity of Algorithm~1 consists of three parts: solving the inverse learning problem, computing the NE, and computing the corresponding Jacobian. Since problem~\eqref{data-driven} is a quadratic program with $M(1+q+r)+N$ decision variables and $MNn$ scalar equality constraints, its complexity is of order
$
\mathcal{O}\!\left((M(Nn+q+r)+N)^3\right).
$ 
Moreover, the computational cost of computing the NE mainly comes from the projection step and is of order $\mathcal{O}(Nn^3)$. The Jacobian computation consists of evaluating the Jacobian of the projection mapping and performing the matrix recursion, with respective complexities $\mathcal{O}(Nn^3)$ and $\mathcal{O}(n^2N^3d)$, yielding an overall cost of $\mathcal{O}(Nn^3+n^2N^3d)$. Therefore, the overall complexity of Algorithm~1 is
$
\mathcal{O}\!\left((M(Nn+q+r)+N)^3+Nn^3+n^2N^3d\right).
$

To steer the NE toward the team-optimal solution, we incorporate an inverse-learning design together with an Adam-type update within a leader--follower framework. Compared with classical gradient-descent-based leader--follower methods under complete information~\cite{grazzi2020iteration,li2023achieving,grontas2024big}, the proposed scheme explicitly addresses a more realistic incomplete-information setting by first learning the unknown follower-side information from observed equilibrium responses. Moreover, compared with gradient descent (GD), which uses only the current first-order hypergradient direction with a common stepsize, the Adam-type update uses first- and second-moment estimates of the hypergradient to adaptively scale the parameter update. This property improves numerical stability, often reduces the number of outer iterations, and hence lowers the total cost of hypergradient evaluations.
In addition, we adopt a non-diminishing stepsize $\eta$ in the Adam update for the mediator's action $\boldsymbol{\theta}$, in contrast to the diminishing-stepsize rule used in~\cite{grontas2024big}. Indeed, diminishing stepsizes may become excessively small as the number of iterations increases, leading to very slow progress in the later stages of the algorithm. By contrast, a non-diminishing stepsize strategy allows the iterates to maintain meaningful updates throughout the computation.


}


\section{Simulation Studies}\label{sec: simulation}

In this section, we illustrate the solution relationships in a shortest path traffic problem \cite{burton1992instance,allen2022using}.
\begin{figure}[ht]
	\centering	
	\includegraphics[scale=0.1]{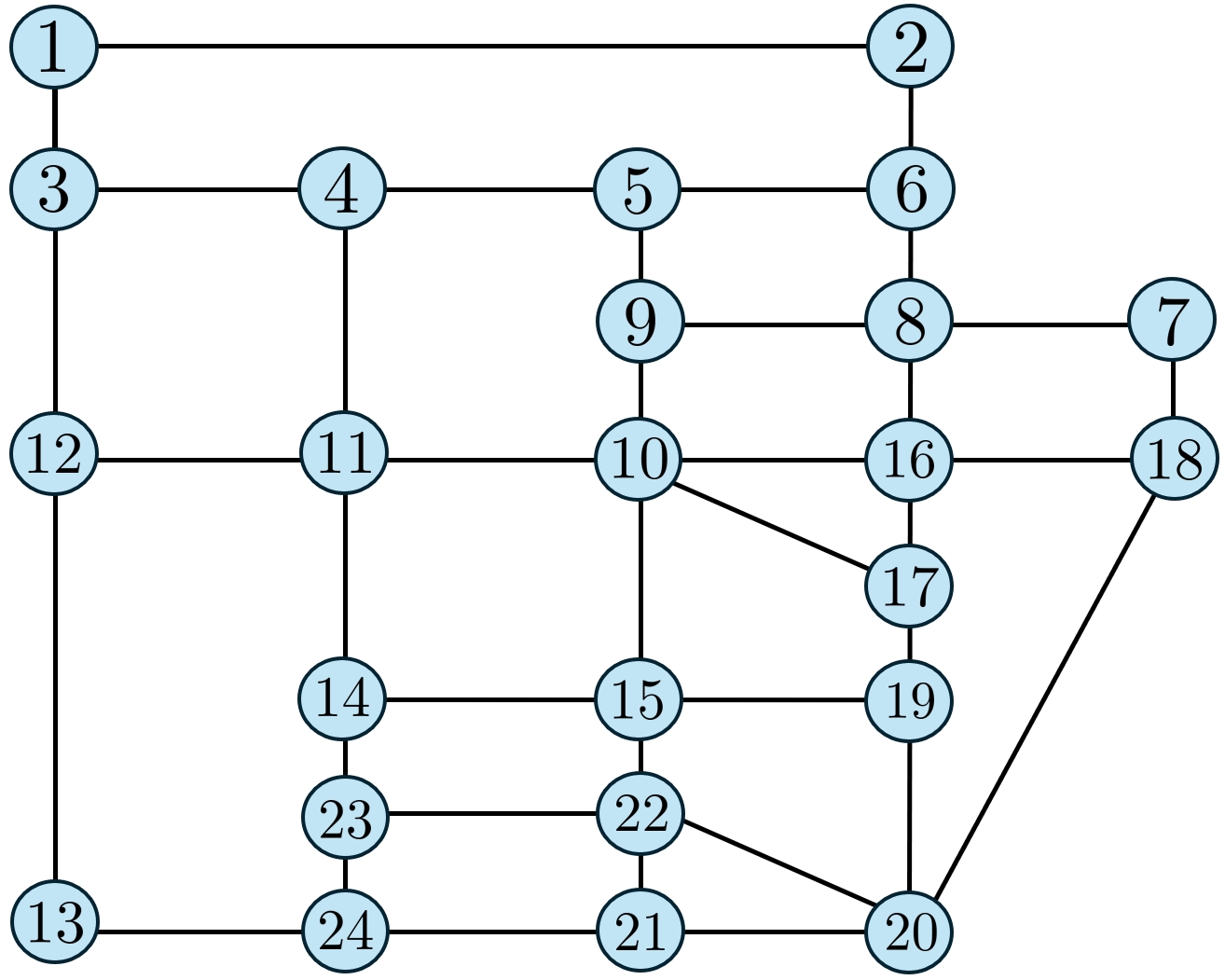}\\
 \vspace{-0.2cm}
	\caption{Sioux Falls network.}
	\label{fig176}
\end{figure}

\noindent\textbf{Model description.} Recall the shortest path traffic problem described in Example 1 in Section \ref{example}.
Consider a network with  $n$ arcs and $r$ nodes, served by $N=4$ vehicles belonging to the same carrier.   Here, we consider the Sioux Falls network \cite{leblanc1975efficient} with $n=31$ arcs and $r=24$ nodes, which is visualized in Fig. \ref{fig176}. For each vehicle $i\in\mathcal{N}=\{1,\dots,4\}$, let $u_i \in \mathbb{R}^n$   be the flow decision variables over the  network arcs.   
Each vehicle’s feasible set is $
{\Xi}_i=\{u_i\in\mathbb{R}^n|u_i\geq \boldsymbol{0}_n, H u_i=m_i\}
    $, where $H\in\{-1,0,1\}^{r\times n}$ is the conservation-of-flow matrix which keeps track of the net flow at the nodes \cite{marcotte2007traffic}, and $m_i\in\mathbb{R}^{r}$
    indicates the origin and destination of vehicle $i$.
    As an example for \(m_i\), if vehicle \(i = 1\) departs from node 3 and arrives at node 18 in the  Sioux Falls network, then \(m_1\) is a \(24\times1\) vector
    with a \(-1\) at entry 3 and a \(+1\) at entry 18,  and zeros elsewhere.
The team objective \eqref{team_traffic} and the individual objective \eqref{game_traffic} are given in Example~1.

\begin{figure}[ht]
\vspace{-0.2cm}
	\centering	
	\includegraphics[scale=0.38]{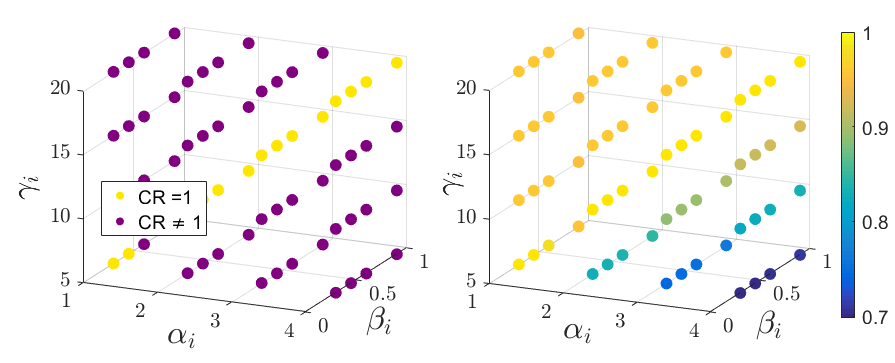}\\
 \vspace{-0.2cm}
	\caption{Left: consistency relation. Right: closeness ratio.
    }
	\label{fig17}
\end{figure}

\noindent\textbf{Experiment results.}  
We begin by examining the consistency and deviation relations, focusing on the case where both the NE and the team-optimal solution are unique.
For ease of illustration, we first consider scalar parameters, i.e. $Q(\alpha)=\alpha I_n$, $Q(\alpha_i)=\alpha_i I_n$  with $\alpha,\alpha_i\in\mathbb{R}$. Similarly, $Q(\beta)=\beta  I_n$,  $Q(\beta_i)=\beta_i I_n$ with $\beta,\beta_i\in\mathbb{R}$, and we take  $\gamma_i*\boldsymbol{1}_n$ with $\gamma_i\in\mathbb{R}$.  
According to \cite{allen2022using}, we fix
$\alpha = 2$, $\beta = 0.3$, $\gamma = 10$,
and then vary each vehicle’s subjective parameter over four values,
$\alpha_i \in \{1,2,3,4\}$,
$\beta_i \in \{0.3,0.45,0.6,0.9\}$,
$\gamma_i \in \{5,10,15,20\}$. For simplicity, we vary all vehicles’ parameters in unison—i.e., we set \(\alpha_i=\alpha_j\), \(\beta_i=\beta_j\), and \(\gamma_i=\gamma_j\) for all \(i,j\in\mathcal{N}\).   
The left subfigure  in Fig.~\ref{fig17} shows the consistency relation (CR) between these two solutions for all \(4\times4\times4=64\) combinations of \(\alpha_i\), \(\beta_i\), and \(\gamma_i\). Here, the CR is indicated by a binary variable.
The yellow points mark cases where $\operatorname{CR} = 1 $ (i.e., consistency), 
while the purple points  mark cases where $\operatorname{CR}\neq 1$ (i.e.,  inconsistency). {The cases with $\operatorname{CR}=1$ correspond to scenarios where the team and individual objectives satisfy the coincidence condition.
} 
Under these consistent cases, each vehicle’s selfish routing strategies also minimizes the total travel time across the entire network, implying self‐interested decision-making causes no additional congestion.

{We further distinguish the $\operatorname{CR}\neq 1$ cases  by considering the degrees of deviation.
We define the closeness ratio as 
$
\operatorname{closeness\;ratio}:=\frac{1}{1+\|\boldsymbol{u}^\Diamond - \boldsymbol{u}^*\|},
$
which ranges from $0$ to $1$.
The right subfigure in Fig.~\ref{fig17} shows a heat map of the closeness ratios under each parameter setting. 
As these parameters vary, the ratio changes over different ranges. The points colored yellow (closeness ratio near 1) denote near‐perfect consistency, indicating smaller deviations between these two solutions, while the points in blue indicate larger deviations. Moreover, Tab.~\ref{tab:bound} compares it with the bound given in Theorem~\ref{t13} for several representative cases. 
Since the bound characterizes a worst-case deviation, it is expected to be larger than the observed deviation. Nevertheless, it remains informative, as larger observed deviations are generally accompanied by larger theoretical bounds.}

\begin{table}[t]
\vspace{-0.2cm}
\centering
\caption{Deviation bound between NE and the team-optimal solution}
\label{tab:bound}
\renewcommand{\arraystretch}{1.08}
\setlength{\tabcolsep}{2.6pt}
\begin{tabular}{c c c c c}
\hline
\textbf{$(\alpha_i,\beta_i,\gamma_i)$} 
& \textbf{$(3,\!0.9,\!20)$} 
& \textbf{$(4,\!0.9,\!15)$} 
& \textbf{$(4,\!0.6,\!10)$} 
& \textbf{$(2,\!0.45,\!5)$} \\
\hline
\textbf{$\|\boldsymbol{u}^{\Diamond}-\boldsymbol{u}^*\|$} 
& 0.17 
& 0.24 
& 0.34 
& 0.58 \\
\textbf{$\xi/{\kappa_1}$} 
& 0.93 
& 1.78 
& 1.68 
& 2.42 \\
\hline
\end{tabular}
\vspace{-0.3cm}
\end{table}

\begin{figure}[ht]
 \vspace{-0.2cm}
	\centering	
	\includegraphics[scale=0.24]{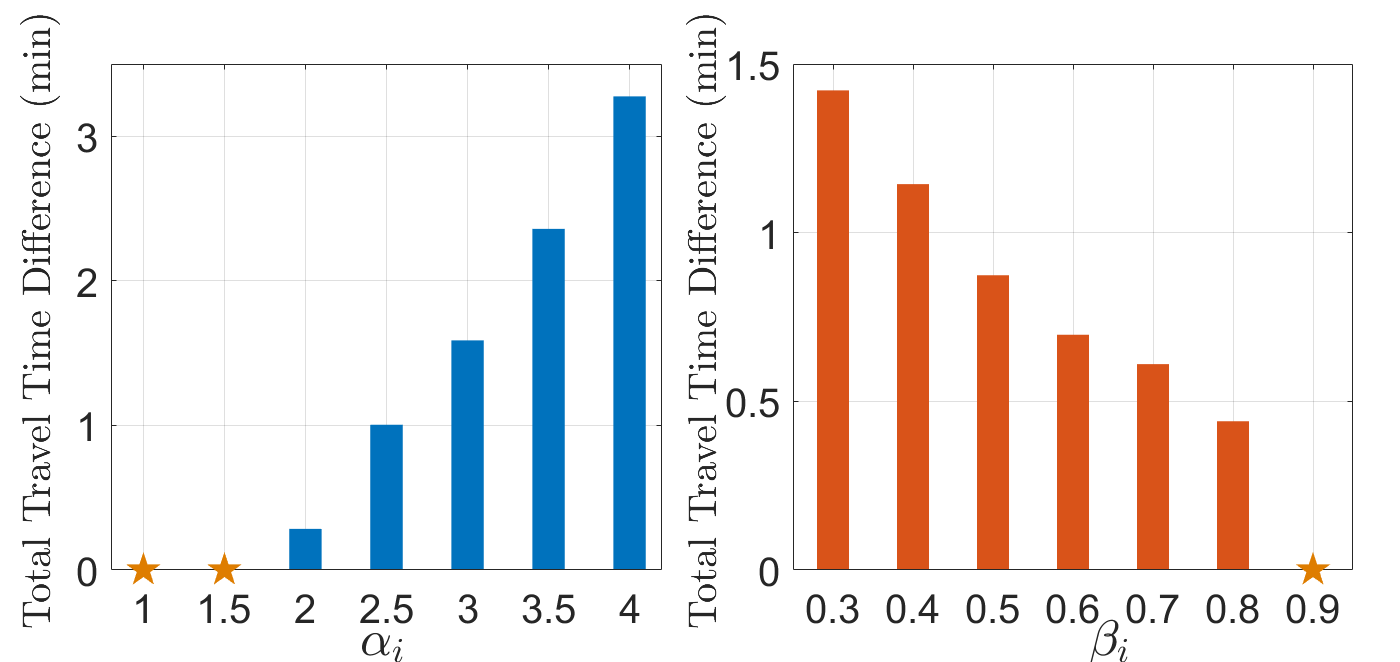}\\
 \vspace{-0.2cm}
	\caption{Travel time difference.}
	\label{fig178}
	\vspace{-0.3cm}
\end{figure}

{Moreover, Fig.~\ref{fig178} shows the impact of different values of $\alpha_i$ and $\beta_i$ on the deviation between the NE and the team-optimal solution, measured by the resulting difference in the total travel time of all vehicles. Recall that $\alpha_i$ weights each vehicle's own-flow cost and  $\beta_i$ weights the interaction with the aggregate flow.
 In the left subfigure, when $\alpha_i$ is small, the consistency relation holds (marked by orange stars), and thus the NE is team-optimal with zero travel time difference. As $\alpha_i$ increases, each vehicle more aggressively avoids arcs that appear costly from its own perspective, which may divert traffic to other already heavily loaded arcs and thereby increase the total network travel time. In the right subfigure, increasing $\beta_i$ reduces the travel time difference. Since $\beta_i$ weights the interaction with the aggregate flow, a larger $\beta_i$ makes each vehicle more responsive to the congestion imposed on and by the rest of the network, thereby promoting a more balanced route allocation and driving the resulting NE closer to the team-optimal solution.}

\begin{figure}[ht]
\vspace{-0.3cm}
\centering	
			\subfigure[ Case under $M=6$]{
				\begin{minipage}[t]{0.49\linewidth}
					\centering
					\includegraphics[width=3.8cm]{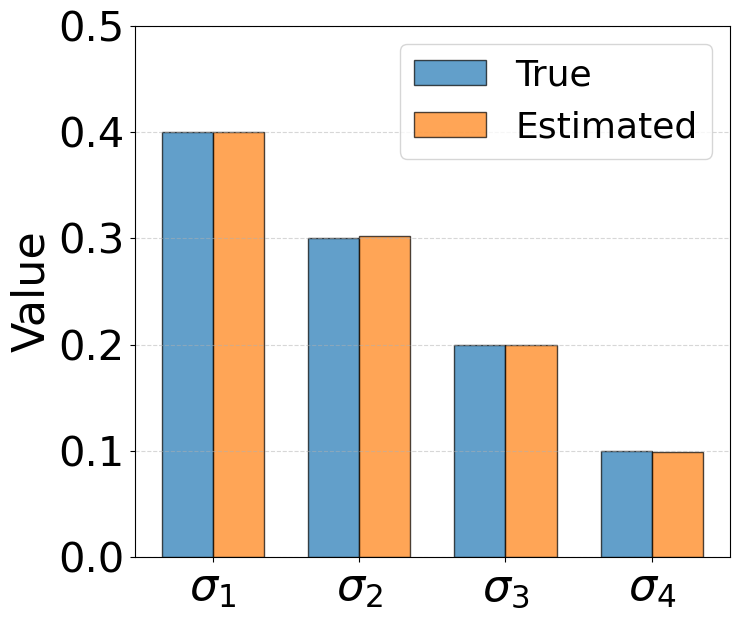}
				\end{minipage}%
			}%
			\subfigure[Different sample cases]{
				\begin{minipage}[t]{0.49\linewidth}
					\centering
					\includegraphics[width=4.5cm]{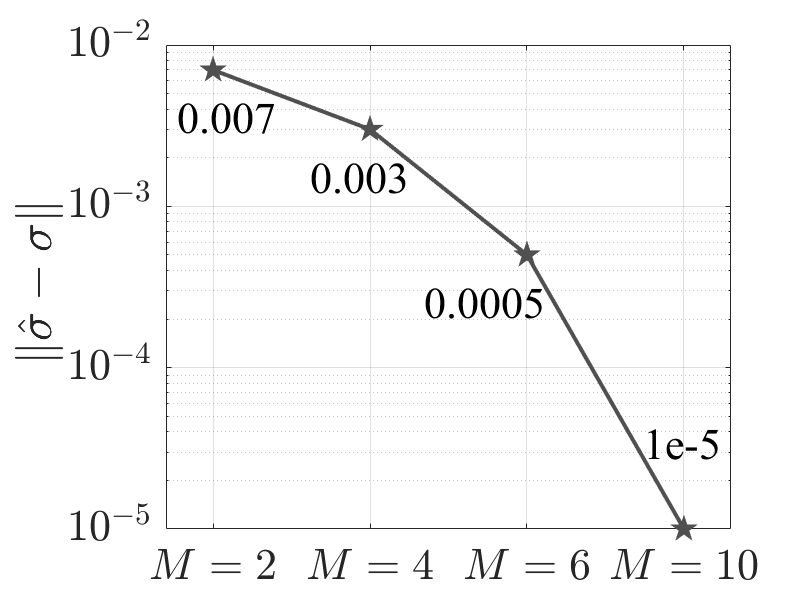}					
				\end{minipage}%
			}\\
 \vspace{-0.2cm}
	\caption{Learning performance}
	\label{fig186}
	\vspace{-0.5cm}
\end{figure}

{ Finally, we verify the effectiveness of the intervention scheme for  shifting NE strategies toward the team‐optimal solution.  Because each vehicle selects its route selfishly based on its own perceived travel cost rather than the team cost, the mediator acts as a network operator that reshapes perceived route costs through interventions to enhance overall system efficiency. We set
$Q(\alpha)\!:=\!\operatorname{diag}(\alpha)$, $Q(\alpha_i)\!:=\!\operatorname{diag}(\alpha_i)$  with $\alpha,\alpha_i\!\in\!\mathbb{R}^n$, $Q(\beta)\!:=\!\operatorname{diag}(\beta)$,  $Q(\beta_i)\!:=\!\operatorname{diag}(\beta_i)$ with $\beta,\beta_i\!\in\!\mathbb{R}^n$, and take  $\gamma_i\!\in \!\mathbb{R}^n$. 
To validate the proposed inverse learning scheme, 
consider that
the true  $\boldsymbol{\sigma}$ is set as
$
\boldsymbol{\sigma}\!=\!(0.4,0.3,0.2,0.1),
$
which is unknown to the mediator.
The mediator first takes different actions and gets the corresponding NE, producing a dataset
$
\{(\boldsymbol{\theta}[j],\boldsymbol{u}^{\Diamond}[j])\}_{j=1}^{M}.
$
where \(M=6\). 
These observed equilibrium samples are then substituted into the proposed inverse problem~\eqref{data-driven} to estimate \(\boldsymbol{\sigma}\).  
Fig.~\ref{fig186}(a) shows that the proposed inverse learning problem is able to recover the unknown weight vector accurately,  where the estimation error is $\|\hat{\boldsymbol{\sigma}}\!-\!\boldsymbol{\sigma}\|=0.0005$. Moreover, Fig.~\ref{fig186}(b) lists the estimation errors under different sample sizes \(M\), showing that as \(M\) increases, the estimation accuracy of \(\boldsymbol{\sigma}\) improves.}

\begin{table*}[t]
\centering
\caption{Ablation study comparing the effects of inverse learning and optimizer choice.}
\label{tab:ablation}
\renewcommand{\arraystretch}{1.15}
\setlength{\tabcolsep}{3.5pt}
\begin{tabular}{l l c c c c}
\hline
\textbf{Method} & \textbf{Model information} & \textbf{Final $\|\boldsymbol{u}^{\Diamond}\!-\!\boldsymbol{u}^*\|$} & \textbf{Final $\mathcal{C}(\boldsymbol{u}^{\Diamond})\!-\!\mathcal{C}(\boldsymbol{u}^*)$} & \textbf{Iterations to tol.} & \textbf{Time to tol.} \\
\hline
No-learning + GD        & misspecified ($\tilde{\boldsymbol{\sigma}}$) & 0.877 & 0.55 & 640 & 1676.75s \\
No-learning + Adam      & misspecified ($\tilde{\boldsymbol{\sigma}}$) & 0.824 & 0.36 & 400 & 1131.55s \\
Inverse-learning + GD   & learned ($\hat{\boldsymbol{\sigma}}$)        & 0.027 & 1e-3 & 636 & 1833.26s \\
\rowcolor{lightgray}
\textbf{Inverse-learning + Adam} & \textbf{learned ($\hat{\boldsymbol{\sigma}}$) }       & \textbf{0.0046} & \textbf{1e-4} & \textbf{340} & \textbf{1219.35s} \\
\hline
\end{tabular}
\end{table*}



\begin{figure}[ht]
   \vspace{-0.2cm}
			\hspace{-0.5cm}
			\centering	
			\subfigure[Case (i)]{
				\begin{minipage}[t]{0.49\linewidth}
					\centering
					\includegraphics[width=4.2cm]{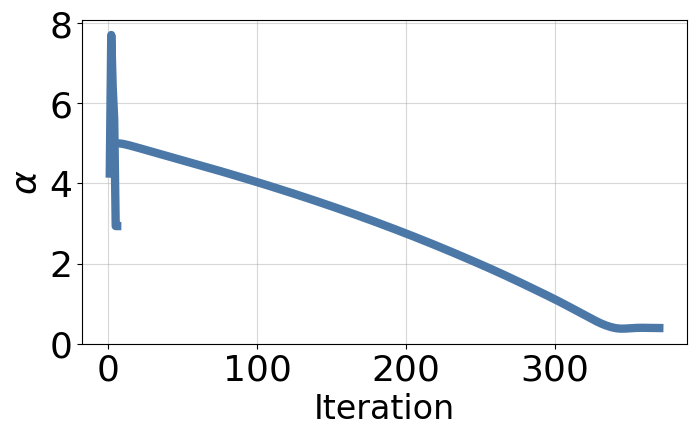}\\
                    \includegraphics[width=4.2cm]{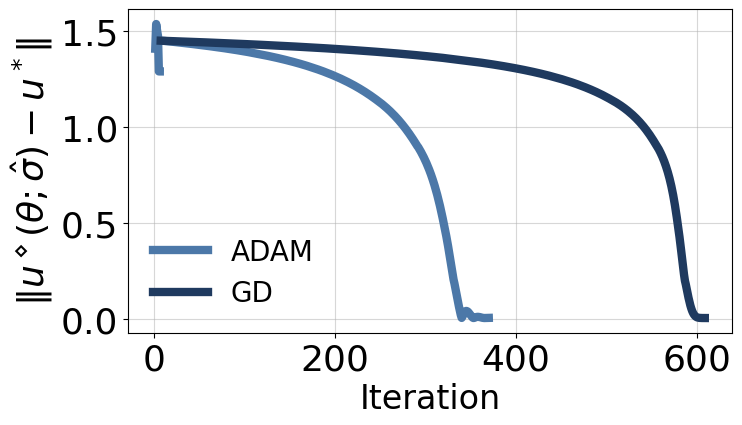}
				\end{minipage}%
			}%
			\subfigure[Case (ii)]{
				\begin{minipage}[t]{0.49\linewidth}
					\centering
					\includegraphics[width=4.2cm]{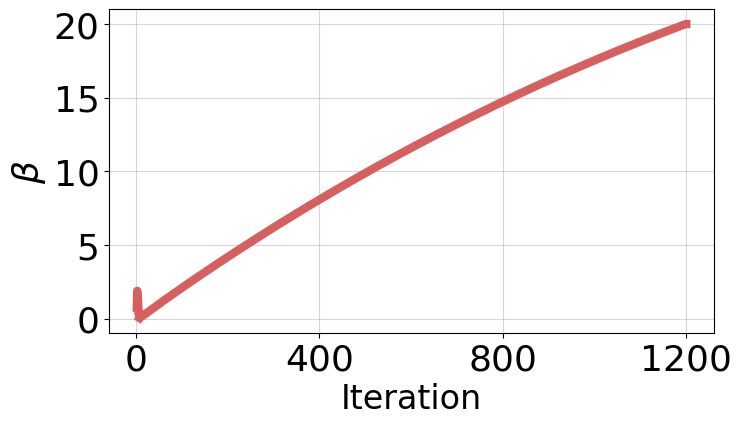}\\
                    \includegraphics[width=4.2cm]{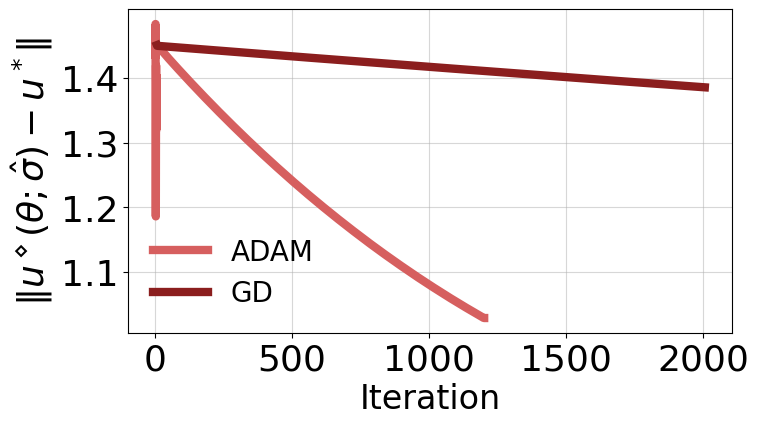}
				\end{minipage}%
			}\\
            \hspace{-0.5cm}
            \subfigure[Case (iii)]{
				\begin{minipage}[t]{0.49\linewidth}
					\centering
					\includegraphics[width=4.2cm]{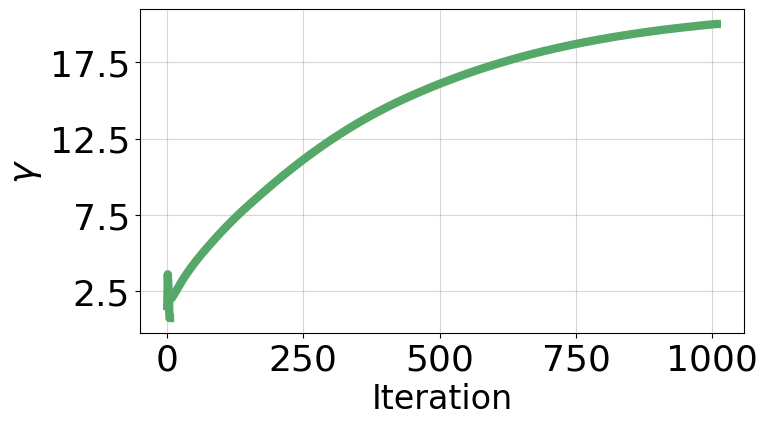}\\
                    \includegraphics[width=4.2cm]{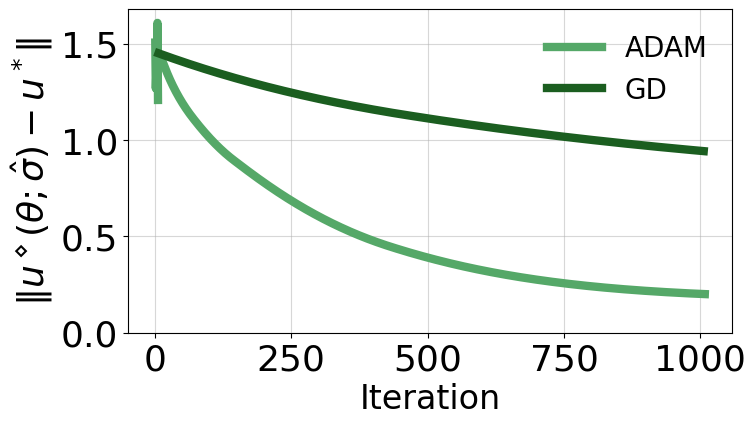}
				\end{minipage}%
			}%
			\subfigure[Case (iv)]{
				\begin{minipage}[t]{0.49\linewidth}
					\centering
					\includegraphics[width=4.2cm]{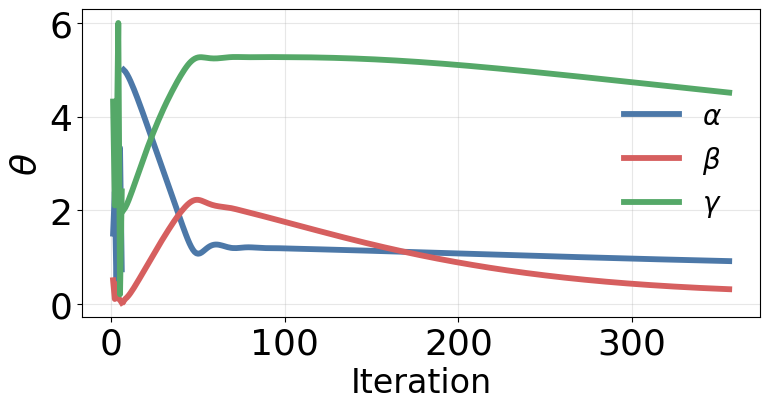}\\
                    \includegraphics[width=4.2cm]{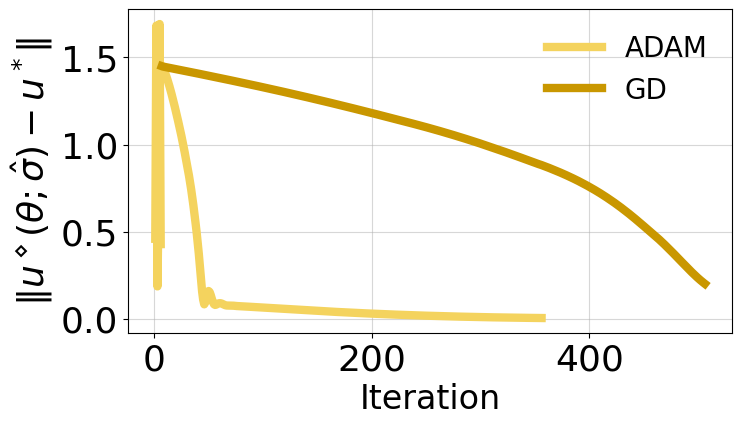}
				\end{minipage}%
			}
			\centering
            \vspace{-0.25cm}
			\caption{Parameter adjustment under different cases.} 
   \vspace{-0.3cm}
			\label{fig_theta}
		\end{figure}

{Then, we validate the Adam scheme. Consider four adjustment scenarios: (1) only  $\{\alpha_i\}_{i=1}^N$ are adjustable with $\boldsymbol{\theta}\!=\!\{\Delta \alpha_i\}_{i=1}^N$;  (2) only the parameters $\{\beta_i\}_{i=1}^N$ are adjustable with $\boldsymbol{\theta}\!=\!\{\Delta \beta_i\}_{i=1}^N$; (3) only the parameters $\{\gamma_i\}_{i=1}^N$ are adjustable with $\boldsymbol{\theta}\!=\!\{\Delta \gamma_i\}_{i=1}^N$; (4) all parameters \(\alpha_i\), \(\beta_i\), and \(\gamma_i\) for \(i\in\mathcal{N}\) are adjustable with $\boldsymbol{\theta}\!=\!\{\Delta \alpha_i,\! \Delta \beta_i,\! \Delta \gamma_i\}_{i=1}^N$.
Fig.~\ref{fig_theta} shows the performance of Alg.~1 under the four adjustment scenarios. In each case, the mediator first estimates $\boldsymbol{\sigma}$ using $M\!=\!6$ samples, and then, starting from iteration $7$, applies the Adam update and the gradient descent (GD) method~\cite{li2023achieving,grontas2024big}, respectively, under the same learned model $\hat{\boldsymbol{\sigma}}$. The stopping criterion is set as $\|\boldsymbol{\theta}[k\!+\!1]\!-\!\boldsymbol{\theta}[k]\|\!\le\! 10^{-5}$. The horizontal axis denotes the iteration number. As shown in Fig.~\ref{fig_theta}, the oscillatory behavior in the first few iterations reflects the random actions taken by the mediator to estimate $\boldsymbol{\sigma}$. After this identification stage, Adam decreases $\|\boldsymbol{u}^{\Diamond}(\boldsymbol{\theta};\hat{\boldsymbol{\sigma}})\!-\!\boldsymbol{u}^*\|$ much faster than standard gradient descent during the subsequent optimization process. Starting from a suitable initial point, Adam reaches a local minimum in substantially fewer outer iterations, whereas GD reduces the error much more gradually.
}

{More specifically, when all parameters are adjustable, the deviation is reduced to a numerically negligible level, indicating that joint adjustment provides the greatest flexibility for aligning the NE with the team-optimal solution. Notably, adjusting only $\alpha_i$ also reduces the deviation to a numerically negligible level, showing that correcting the vehicles' excessive emphasis on their own local cost can by itself recover the team-optimal solution. Adjusting only $\gamma_i$  yields a very small residual gap and therefore also provides an effective intervention mechanism. By contrast, adjusting only $\beta_i$ still reduces the deviation, but its effect is comparatively weaker. 

Finally, we conduct an ablation study to demonstrate the effectiveness of Alg.~1 from both the information and optimization perspectives. Specifically, we compare four variants: GD and Adam under a misspecified aggregation model without inverse learning, and GD and Adam under the learned aggregation model obtained from~\eqref{data-driven}. For the no-learning baselines, the mediator directly performs the upper-level update based on a misspecified aggregation vector $\tilde{\boldsymbol{\sigma}}=\frac{1}{N}\mathbf{1}_N$, without estimating the true $\boldsymbol{\sigma}$. For the inverse-learning variants, the mediator first estimates $\hat{\boldsymbol{\sigma}}$ from the observed equilibrium samples and then performs the upper-level parameter adjustment. To ensure a fair comparison, all four variants use the same initialization, stopping criterion, and lower-level equilibrium solver. Tab.~\ref{tab:ablation} reports the final steering error $\|\boldsymbol{u}^{\Diamond}-\boldsymbol{u}^*\|$, the final team-cost gap $\mathcal{C}(\boldsymbol{u}^{\Diamond})-\mathcal{C}(\boldsymbol{u}^*)$, the number of outer iterations required to reach a prescribed tolerance $\|\boldsymbol{\theta}[k+1]-\boldsymbol{\theta}[k]\| \leq 1e-5$, and the  computation time (sec). The comparison between the no-learning and inverse-learning variants illustrates the value of estimating the unknown aggregation weights, while the comparison between GD and Adam under the same information pattern isolates the effect of the optimizer itself. As shown in Tab.~\ref{tab:ablation}, the inverse-learning variants consistently outperform their no-learning counterparts, confirming that a behaviorally consistent game model is important for effective upper-level intervention under incomplete information.  Moreover, under both the misspecified and learned models, Adam requires fewer outer iterations and less computation time than GD to achieve the same accuracy level, indicating its superior numerical efficiency. Among all four methods, the inverse-learning + Adam variant achieves the smallest final steering error and team-cost gap, which confirms the benefit of combining model learning with adaptive hypergradient-based optimization.
}



\section{Conclusions}\label{sec: conclusion}
In this paper, we studied the relationship between team-optimal solutions and NE in a class of static team and game problems, aiming to assess the impact of selfish behavior on team performance. We established necessary and sufficient conditions for every NE to be team-optimal and derived an upper bound on the deviation when the consistency does not hold. We further proposed a leader-follower framework that combines inverse learning with an Adam-type update to steer NE strategies toward team optimality. Future work will extend the analysis to dynamic information structures, multi-stage decision-making, and stochastic environments.



\begin{appendices}
\setcounter{section}{0}

\setcounter{mythm}{0}   
\renewcommand{\themythm}{S\arabic{mythm}}
\setcounter{assumption}{0}   
\renewcommand{\theassumption}{C\arabic{assumption}}
\setcounter{lemma}{0}   
\renewcommand{\thelemma}{C\arabic{lemma}}
\renewcommand{\theassumption}{C\arabic{assumption}}

\section{Proof of Lemma \ref{l11}}\label{jaco_con}
Recall that \( \!\phi(\varsigma(\boldsymbol{\theta}, \boldsymbol{u}^\Diamond(\boldsymbol{\theta};\hat{\boldsymbol{\sigma}})))=\zeta(\boldsymbol{\theta},\boldsymbol{u}^{\Diamond}(\boldsymbol{\theta};\hat{\boldsymbol{\sigma}}))\).  Since \(\phi(\cdot)\) is continuously differentiable at
\(\varsigma(\boldsymbol{\theta},\boldsymbol{u}^\Diamond(\boldsymbol{\theta};\hat{\boldsymbol{\sigma}}))\),
and \(F\) is continuously differentiable, it follows that \( \zeta(\boldsymbol{\theta}, \boldsymbol{u}^\Diamond(\boldsymbol{\theta};\hat{\boldsymbol{\sigma}})) \) is continuously differentiable at \( \boldsymbol{u}^\Diamond(\boldsymbol{\theta};\hat{\boldsymbol{\sigma}}) \). 
Moreover, by Lemma \ref{con_ne}, the mapping \( \zeta(\boldsymbol{\theta}, \cdot) \) is  Lipschitz continuous for any \( \boldsymbol{\theta} \), and hence, by~\cite[Theorem 2.17(i)]{outrata2013nonsmooth}, we have
$
\left\| \text{J}_{\boldsymbol{u}^{\Diamond}}{\zeta}(\boldsymbol{\theta}, \boldsymbol{u}^\Diamond (\boldsymbol{\theta};\hat{\boldsymbol{\sigma}})) \right\| < 1.
$ 
Let us define the mapping $\Phi(\boldsymbol{\theta}, \boldsymbol{u};\hat{\boldsymbol{\sigma}}):\mathbb{R}^{dN}\times \mathbb{R}^{nN}\rightarrow \mathbb{R}^{nN}$,
\begin{equation}\label{phi}
\Phi(\boldsymbol{\theta}, \boldsymbol{u}) := \boldsymbol{u} - \zeta(\boldsymbol{\theta}, \boldsymbol{u};\hat{\boldsymbol{\sigma}}),
\end{equation}
such that \( \Phi(\boldsymbol{\theta}, \boldsymbol{u}^\Diamond(\boldsymbol{\theta};\hat{\boldsymbol{\sigma}})) = 0 \).
According to the implicit function theorem~\cite[Theorem 1B.1]{dontchev2009implicit}, since $\Phi$ is continuously differentiable at \( (\boldsymbol{\theta},\boldsymbol{u}^\Diamond(\boldsymbol{\theta};\hat{\boldsymbol{\sigma}})) \), 
and \(  I - \text{J}_{\boldsymbol{u}^{\Diamond}}{\zeta}(\boldsymbol{\theta}, \boldsymbol{u}^\Diamond (\boldsymbol{\theta};\hat{\boldsymbol{\sigma}})) \) is nonsingular, it follows that \( \boldsymbol{u}^\Diamond(\boldsymbol{\theta};\hat{\boldsymbol{\sigma}}) \) is continuously differentiable at \( \boldsymbol{\theta} \), with Jacobian given by $\text{J}{\boldsymbol{u}^\Diamond}(\boldsymbol{\theta};\hat{\boldsymbol{\sigma}})\! =\! \!\left( I \!\!-\! \text{J}_{\boldsymbol{u}^{\Diamond}}{\zeta}(\boldsymbol{\theta}, \!\boldsymbol{u}^\Diamond(\boldsymbol{\theta};\!\hat{\boldsymbol{\sigma}})) \right)^{-1}\! \!\text{J}_{\boldsymbol{\theta}}{\zeta}(\boldsymbol{\theta}, \!\boldsymbol{u}^\Diamond (\boldsymbol{\theta};\!\hat{\boldsymbol{\sigma}})).$
Furthermore, because $
\left\| \text{J}_{\boldsymbol{u}^{\Diamond}}{\zeta}(\boldsymbol{\theta}, \boldsymbol{u}^\Diamond (\boldsymbol{\theta};\!\hat{\boldsymbol{\sigma}})) \right\| < 1
$, the update in \eqref{jaco} is a contraction and therefore converges to a unique fixed point, which satisfies the expression of $\text{J}{\boldsymbol{u}^\Diamond}(\boldsymbol{\theta};\hat{\boldsymbol{\sigma}})$.
\hfill $\square$

{\section{Proof of Theorem \ref{t5}}\label{proofthm3}

To establish the convergence of~\eqref{adam}, we first show that the function \( \boldsymbol{u}^\Diamond(\boldsymbol{\theta};\!\hat{\boldsymbol{\sigma}}) \) is definable. Since definability is preserved under function composition~\cite{bolte2021nonsmooth}, this also implies that both \( \psi \) and
$\Psi(\boldsymbol{\theta};\hat{\boldsymbol{\sigma}})
:=
\psi(\boldsymbol{\theta};\!\hat{\boldsymbol{\sigma}})+\frac{\rho}{2}\|\boldsymbol{\theta}\|^2$
are definable. This allows us to invoke the differential inclusion framework from~\cite[Theorem 3.2]{davis2020stochastic}.
Recall from~\eqref{phi} that \( \boldsymbol{u}^\Diamond(\boldsymbol{\theta};\!\hat{\boldsymbol{\sigma}}) \) is implicitly defined by the mapping \( \Phi(\boldsymbol{\theta},\boldsymbol{u};\!\hat{\boldsymbol{\sigma}}) \). Note that the projection operator \( \phi(\cdot) \) is piecewise affine and thus semi-algebraic. Moreover, the mapping \(F\) is definable by Assumption~\ref{assum3}. Hence,
$\zeta(\boldsymbol{\theta},\boldsymbol{u};\!\hat{\boldsymbol{\sigma}})
=
\phi\bigl(\boldsymbol{u}-\tau F(\boldsymbol{\theta},\boldsymbol{u};\hat{\boldsymbol{\sigma}})\bigr)$
is definable, and therefore \( \Phi \) is definable as well. By the non-smooth implicit function theorem~\cite[Theorem 5]{bolte2021nonsmooth}, it follows that \( \boldsymbol{u}^\Diamond(\boldsymbol{\theta};\hat{\boldsymbol{\sigma}}) \) is definable. Consequently, \(\psi\) and \(\Psi\) are definable. Furthermore, by the conservative Jacobian chain rule~\cite[Proposition 2]{bolte2021nonsmooth}, the conservative gradient of \(\psi\) is given by
$\mathcal{J}_\psi(\boldsymbol{\theta};\hat{\boldsymbol{\sigma}})
=
\{
\tilde{J}^{\top}\bigl(\boldsymbol{u}^{\Diamond}(\boldsymbol{\theta};\hat{\boldsymbol{\sigma}})-\boldsymbol{u}^{*}\bigr)
\;:\;
\tilde{J}\in\mathcal{J}_{\boldsymbol{u}^{\Diamond}}(\boldsymbol{\theta};\hat{\boldsymbol{\sigma}})
\}$,
where \(\mathcal{J}_{\boldsymbol{u}^{\Diamond}}(\boldsymbol{\theta};\hat{\boldsymbol{\sigma}})\) denotes the conservative Jacobian of \(\boldsymbol{u}^{\Diamond}\) at \(\boldsymbol{\theta}\).

We now apply the differential inclusion framework to analyze the convergence of the iterates \( \{\boldsymbol{\theta}[k]\} \).
To this end, introduce the auxiliary sequence
$y[k]:=-\frac{1}{\rho}\mu[k], \forall k\ge 0.$
Substituting it into the update of \(\mu[k]\) in~\eqref{adam}, we obtain
\begin{equation}\label{inter1}
y[k+1]
=
y[k]-\frac{\vartheta_1[k]}{\rho}\bigl(\omega[k]+\rho y[k]\bigr).
\end{equation}
Consider the differential inclusion
\begin{equation}\label{inclusion_y2}
\frac{dy}{dt}\in -\mathcal{J}_{\Psi}(y;\hat{\boldsymbol{\sigma}}).
\end{equation}

We then use the following consequence of \cite[Proposition 3.27, Theorem 3.6]{benaim2005stochastic} and \cite[Theorem 3.2]{davis2020stochastic}.

\begin{lemma}\label{ly2}
Under Assumption~\ref{assum3}, let \( \{y[k]\}_{k\in\mathbb N} \) be generated by~\eqref{inter1}, and suppose that
\begin{enumerate}
    \item the set \(\{\Psi(y):y\in \operatorname{crit}_{\Psi}\}\) has empty interior in \(\mathbb R\);
    \item \(\Psi\) is a Lyapunov function for the differential inclusion~\eqref{inclusion_y2}, with stable set \(\operatorname{crit}_{\Psi}\);
    \item the linearly interpolated process associated with \(\{y[k]\}_{k\in\mathbb N}\) is a perturbed solution of~\eqref{inclusion_y2}.
\end{enumerate}
Then every limit point of \(\{y[k]\}\) lies in \(\operatorname{crit}_{\Psi}\), and the sequence \(\{\Psi(y[k];\hat{\boldsymbol{\sigma}})\}\) converges.
\end{lemma}

It remains to verify the three conditions in Lemma~\ref{ly2}. First, Condition~(1) is a direct consequence of the non-smooth weak Sard property~\cite{bolte2021conservative,davis2020stochastic}. Since \(\Psi\) and \(\mathcal{J}_{\Psi}\) are definable, \cite[Theorem 5]{bolte2021conservative} implies that \(\Psi\) is path-differentiable and that the set \(\{\Psi(y;\hat{\boldsymbol{\sigma}}):y\in\operatorname{crit}_{\Psi}\}\) is finite. Hence Condition~(1) holds.
Second, Condition~(2) follows from \cite{bolte2021conservative}, which guarantees that a definable path-differentiable function is a Lyapunov function for the corresponding conservative differential inclusion.  Finally, we verify Condition~(3). Since \(\omega[k]\in \mathcal J_\psi(\boldsymbol{\theta}[k];\hat{\boldsymbol{\sigma}})\), we have \(\omega[k]+\rho\boldsymbol{\theta}[k]\in \mathcal J_\Psi(\boldsymbol{\theta}[k];\hat{\boldsymbol{\sigma}})\). Moreover,
$-(\omega[k]+\rho y[k])=-(\omega[k]+\rho\boldsymbol{\theta}[k])-\rho\bigl(y[k]-\boldsymbol{\theta}[k]\bigr)$,
where the first term belongs to \(-\mathcal J_\Psi(\boldsymbol{\theta}[k];\hat{\boldsymbol{\sigma}})\). Thus, it suffices to show that \(\rho(y[k]-\boldsymbol{\theta}[k])\to 0\). To this end, define \(\Delta[k]:=\rho\boldsymbol{\theta}[k]+\mu[k]\). Using the updates in~\eqref{adam}, we obtain \(\Delta[k+1]=\bigl(1-\eta\rho(\sqrt{s[k+1]}+\varpi)^{-1}\bigr)\odot\bigl(\Delta[k]+\vartheta_1[k](\omega[k]-\mu[k])\bigr)\). Since \(0<\eta<\varpi/\rho\), there exists \(\tilde\eta\in(0,1)\) such that \(\|1-\eta\rho(\sqrt{s[k+1]}+\varpi)^{-1}\|_\infty\le 1-\tilde\eta\) for all \(k\ge 0\). Moreover, by Assumption~\ref{assum3}(4), together with the Adam recursions, the sequences \(\{\omega[k]\}\), \(\{\mu[k]\}\), \(\{s[k]\}\), and \(\{\boldsymbol{\theta}[k]\}\) are bounded. Hence there exists \(M_d>0\) such that \(\sup_{k\ge 0}\|\omega[k]-\mu[k]\|\le M_d\). Therefore,
$\|\Delta[k+1]\|\le (1-\tilde\eta)\|\Delta[k]\|+M_d\vartheta_1[k]$.
Iterating this inequality gives
$\|\Delta[k+1]\|\le (1-\tilde\eta)^{k+1}\|\Delta[0]\|+M_d\sum_{i=0}^k(1-\tilde\eta)^{k-i}\vartheta_1[i]$.
Since \(\vartheta_1[k]\to 0\), the convolution term converges to \(0\), and thus \(\Delta[k]\to 0\). Recalling that \(y[k]=-\mu[k]/\rho\), we obtain \(\|\boldsymbol{\theta}[k]-y[k]\|=\rho^{-1}\|\Delta[k]\|\to 0\), and hence \(\rho(y[k]-\boldsymbol{\theta}[k])\to 0\). Since \(\mathcal J_\Psi\) is outer semicontinuous and locally bounded, it follows that the linearly interpolated process associated with \(\{y[k]\}\) is a perturbed solution of the differential inclusion~\eqref{inclusion_y2}. Hence Condition~(3) of Lemma~\ref{ly2} holds. Applying Lemma~\ref{ly2}, we conclude that every limit point of \(\{y[k]\}\) lies in \(\operatorname{crit}_\Psi\), and \(\{\Psi(y[k];\hat{\boldsymbol{\sigma}})\}\) converges.

Finally, since \(\|\boldsymbol{\theta}[k]-y[k]\|\to 0\), every limit point of \(\{\boldsymbol{\theta}[k]\}\) coincides with a limit point of \(\{y[k]\}\). Therefore every limit point of \(\{\boldsymbol{\theta}[k]\}\) lies in \(\operatorname{crit}_{\Psi}\). Moreover, since \(\Psi\) is continuous and \(\|\boldsymbol{\theta}[k]-y[k]\|\to 0\), we have
$|\Psi(\boldsymbol{\theta}[k];\hat{\boldsymbol{\sigma}})-\Psi(y[k];\hat{\boldsymbol{\sigma}})|\to 0$.
Because \(\{\Psi(y[k];\hat{\boldsymbol{\sigma}})\}\) converges, it follows that \(\{\Psi(\boldsymbol{\theta}[k];\hat{\boldsymbol{\sigma}})\}\) also converges. 
\hfill $\square$}

\section{Proof of Corollary \ref{cor:const_momentum}}
{
As in the proof of Theorem~\ref{t5}, the lower-level solution mapping \(\boldsymbol{u}^\Diamond(\boldsymbol{\theta};\hat{\boldsymbol{\sigma}})\) is definable. Hence \(\psi\) and \(\Psi(\boldsymbol{\theta};\hat{\boldsymbol{\sigma}})\) are definable. By~\cite[Theorem~5]{bolte2021conservative}, \(\Psi\) is path-differentiable, satisfies the weak Sard property, and is a Lyapunov function for the differential inclusion \(\dot y\in-\mathcal J_\Psi(y;\hat{\boldsymbol{\sigma}})\).

Define \(y[k]:=-\mu[k]/\rho\). Since \(\vartheta_1[k]\equiv\bar\vartheta\), the update of \(\mu[k]\) in~\eqref{adam} becomes \(y[k+1]=y[k]-\frac{\bar\vartheta}{\rho}\bigl(\omega[k]+\rho y[k]\bigr)\), where \(\omega[k]\in\mathcal J_\psi(\boldsymbol{\theta}[k];\hat{\boldsymbol{\sigma}})\). Thus, \(\{y[k]\}\) is a constant-stepsize discretization of \(\dot y\in-\mathcal J_\Psi(y;\hat{\boldsymbol{\sigma}})\), up to the discrepancy between \(\boldsymbol{\theta}[k]\) and \(y[k]\). Let \(e[k]:=\boldsymbol{\theta}[k]-y[k]\). Since \(\mu[k+1]=-\rho y[k+1]\), the update of \(\boldsymbol{\theta}[k]\) yields \(e[k+1]=\bigl(1-\eta\rho(\sqrt{s[k+1]}+\varpi)^{-1}\bigr)\odot\bigl(e[k]+y[k]-y[k+1]\bigr)\). Since \(0<\eta<\varpi/\rho\) and \(\{s[k]\}\) is bounded, there exists \(q\in(0,1)\) such that \(\bigl\|1-\eta\rho(\sqrt{s[k+1]}+\varpi)^{-1}\bigr\|_\infty\le q\) for all \(k\). Moreover, \(\{y[k]\}\) and \(\{\omega[k]\}\) are bounded, so \(\|y[k+1]-y[k]\|=\frac{\bar\vartheta}{\rho}\|\omega[k]+\rho y[k]\|\le C\bar\vartheta\) for some \(C>0\). Hence \(\|e[k+1]\|\le q\|e[k]\|+qC\bar\vartheta\), which implies \(\limsup_{k\to\infty}\|e[k]\|\le \frac{qC}{1-q}\bar\vartheta\). Therefore, for any \(\delta>0\), there exists \(\hat\vartheta_0(\delta)\in(0,1)\) such that \(0<\bar\vartheta\le \hat\vartheta_0(\delta)\) implies \(\limsup_{k\to\infty}\|\boldsymbol{\theta}[k]-y[k]\|\le \delta\).

Next, since \(\omega[k]+\rho\boldsymbol{\theta}[k]\in\mathcal J_\Psi(\boldsymbol{\theta}[k];\hat{\boldsymbol{\sigma}})\), the recursion of \(y[k]\) can be written as \(y[k+1]=y[k]-\frac{\bar\vartheta}{\rho}\bigl(v[k]+r[k]\bigr)\), where \(v[k]\in\mathcal J_\Psi(y[k];\hat{\boldsymbol{\sigma}})\) and \(r[k]\) is the associated error term. By the outer semicontinuity and local boundedness of \(\mathcal J_\Psi\), the bound \(\limsup_{k\to\infty}\|\boldsymbol{\theta}[k]-y[k]\|\le \delta\) implies \(\limsup_{k\to\infty}\|r[k]\|\le \nu\) for arbitrarily small \(\nu>0\), provided that \(\bar\vartheta\) is sufficiently small. Therefore, \(\{y[k]\}\) is a small-stepsize, small-perturbation discretization of \(\dot y\in-\mathcal J_\Psi(y;\hat{\boldsymbol{\sigma}})\). By~\cite[Theorem~3.13]{ding2024adamfamily}, for any \(\varepsilon>0\), there exists \(\hat\vartheta_1(\varepsilon/2)\in(0,1)\) such that \(0<\bar\vartheta\le \hat\vartheta_1(\varepsilon/2)\) implies \(\limsup_{k\to\infty}\operatorname{dist}(y[k],\operatorname{crit}_\Psi)\le \varepsilon/2\). Let \(\hat\vartheta_2(\varepsilon/2):=\hat\vartheta_0(\varepsilon/2)\) and define \(\hat\vartheta(\varepsilon):=\min\{\hat\vartheta_1(\varepsilon/2),\hat\vartheta_2(\varepsilon/2)\}\). Then, whenever \(0<\bar\vartheta\le \hat\vartheta(\varepsilon)\), we have both \(\limsup_{k\to\infty}\|\boldsymbol{\theta}[k]-y[k]\|\le \varepsilon/2\) and \(\limsup_{k\to\infty}\operatorname{dist}(y[k],\operatorname{crit}_\Psi)\le \varepsilon/2\). The conclusion follows from \(\operatorname{dist}(\boldsymbol{\theta}[k],\operatorname{crit}_\Psi)\le \|\boldsymbol{\theta}[k]-y[k]\|+\operatorname{dist}(y[k],\operatorname{crit}_\Psi)\). \hfill\(\square\)
}

\end{appendices}

\section*{References}
	\bibliographystyle{IEEEtran}
\bibliography{reference,IDS}

\end{document}